\documentclass[12pt,reqno]{article}
\usepackage{microtype}
\usepackage{colordvi}
\usepackage{extpfeil}
\usepackage{amsmath,amsfonts,amsthm,amssymb,latexsym}
\usepackage{indentfirst}
\allowdisplaybreaks[4]
\usepackage{graphicx}
\usepackage{verbatim}
\usepackage{color}
\usepackage{amscd}
\usepackage{verbatim}
\usepackage{graphicx}
\usepackage{cleveref}

\usepackage{url}

\usepackage{setspace}

\usepackage[left=1.25in, right=1.25in, top=1in, bottom=1in]{geometry}
\usepackage{longtable}

\usepackage{CJKutf8}
\usepackage{appendix}
\usepackage{amsmath,amssymb,amsthm,color,mathrsfs}
\usepackage{enumitem,anysize}%
\usepackage{verbatim}
\usepackage{amssymb}
\usepackage{tikz}
\usepackage{pgfplots}
\marginsize{1in}{1in}{1in}{1in}%

\newtheorem{theorem}{Theorem}[section]
\newtheorem{corollary}[theorem]{Corollary}

\newtheorem{Thm}[theorem]{Theorem}

\newtheorem{lemma}[theorem]{Lemma}

\newtheorem{example}[theorem]{Example}

\newtheorem{problem}[theorem]{Problem}
\newtheorem{definition}[theorem]{Definition}
\newtheorem{Def}[theorem]{Definition}

\newtheorem{assumption}[theorem]{Assumption}
\newtheorem{assumptions}[theorem]{Assumptions}

\newtheorem{remark}[theorem]{Remark}

\newtheorem{examples}[theorem]{Examples}
\newtheorem{remarks}[theorem]{Remarks}
\newtheorem{claim*}[theorem]{Claim}
\newtheorem{notation}[theorem]{Notation}

\newcommand{\sms}{\setminus}
\newcommand{\R}{\mathbb{R}}

\newcommand{\N}{\mathbb{N}}

\newcommand{\blemma}{\begin{lemma}}
	\newcommand{\elemma}{\end{lemma}}
\newcommand{\bnotation}{\begin{notation}}
	\newcommand{\enotation}{\end{notation}}
\newcommand{\bproof}{\begin{proof}}
	\newcommand{\eproof}{\end{proof}}
\newcommand{\bremark}{\begin{remark}}
	\newcommand{\eremark}{\end{remark}}
    \newcommand{\bremarks}{\begin{remarks}}
	\newcommand{\eremarks}{\end{remarks}}
\newcommand{\bcorollary}{\begin{corollary}}
	\newcommand{\ecorollary}{\end{corollary}}
\newcommand{\btheorem}{\begin{theorem}}
	\newcommand{\etheorem}{\end{theorem}}
 \newcommand{\bproblem}{\begin{problem}}
	\newcommand{\eproblem}{\end{problem}}
\newcommand{\bdefinition}{\begin{definition}}
	\newcommand{\edefinition}{\end{definition}}
\newcommand{\bequation}{\begin{equation}}
	\newcommand{\eequation}{\end{equation}}
\newcommand{\bequationn}{\begin{equation*}}
	\newcommand{\eequationn}{\end{equation*}}
\newcommand{\beqnarray}{\begin{eqnarray}}
	\newcommand{\eeqnarray}{\end{eqnarray}} \newcommand{\beqnarrayn}{\begin{eqnarray*}}
	\newcommand{\eeqnarrayn}{\end{eqnarray*}}
 \newcommand{\bexample}{\begin{example}}
	\newcommand{\eexample}{\end{example}}
     \newcommand{\bexamples}{\begin{examples}}
	\newcommand{\eexamples}{\end{examples}}
  \newcommand{\bassumption}{\begin{assumption}}
	\newcommand{\eassumption}{\end{assumption}}
\newcommand{\loc}{{\rm loc}}
\numberwithin{equation}{section}

\newcommand\blue[1]{\textcolor{blue}{#1}}

\def\bal#1\eal{\begin{align}#1\end{align}}
\def\bals#1\eals{\begin{align*}#1\end{align*}}

\newcommand{\Hmm}[1]{\leavevmode{\marginpar{\tiny%
			$\hbox to 0mm{\hspace*{-0.5mm}$\leftarrow$\hss}%
			\vcenter{\vrule depth 0.1mm height 0.1mm width \the\marginparwidth}%
			\hbox to
			0mm{\hss$\rightarrow$\hspace*{-0.5mm}}$\\\relax\raggedright #1}}}

\DeclareMathOperator{\diam}{diam}
\DeclareMathOperator{\supp}{supp}
\DeclareMathOperator{\capacity}{Cap}
\DeclareMathOperator{\vol}{Vol}
\DeclareMathOperator{\sgn}{sgn}
\newcommand{\core}{C_c^{\infty}(\Omega)}

\newcommand{\dHnn}{\,\mathrm{d}\mathcal{H}^{n-1}}
\newcommand{\dx}{\,\mathrm{d}x}

\newcommand{\dt}{\,\mathrm{d}t}

\newcommand{\myd}{\displaystyle}

\DeclareMathOperator{\dive}{div}

\def\<{\langle}
\def\>{\rangle}

\long\def\prob#1\soln#2\endps{{\color{blue}#1}\medskip\par
	\noindent\underline{\sc Solution}:\hspace*{1em}\parindent=2em #2}

\AtBeginDocument{\begin{CJK}{UTF8}{gbsn}}
	\AtEndDocument{\end{CJK}}

      \def\gw{\omega}
                \def\gz{\zeta}
\def\Gg{\Gamma}

\def\Gw{\Omega}              
\pgfplotsset{compat=1.18}
\begin{document}
	\pagenumbering{gobble}
	\title{\textbf{Fuchsian-type singularity for the Finsler $p$-Laplacian with potential}}
\author{Yongjun Hou\thanks{Department of Mathematics, Technion -
Israel Institute of Technology, Haifa 3200003, Israel; yongjun.hou@campus.technion.ac.il; houmathlaw@outlook.com}}
 \date{June 23, 2026}
	\maketitle
	\pagenumbering{arabic}
	\vspace{-10mm}	
 	\begin{abstract}
Let $\Omega\subseteq\mathbb{R}^{n}$ ($2\leq n\in\mathbb{N}$) be a domain and let $\zeta\in\{0,\infty\}$ be an isolated point of the boundary of $\Omega$ in the one-point compactification of $\mathbb{R}^{n}$ with the ideal point $\infty$. Under some further conditions, we study Fuchsian-type singularity at $\zeta$ for the Finsler $p$-Laplace equation with a potential 
 $$-\mathrm{div}\mathcal{A}(x,\nabla u)+V|u|^{p-2}u=0\quad (1<p<\infty)\qquad \mbox{in } \Omega,$$
where $\mathcal{A}(x,\xi)\triangleq\nabla_{\xi}(H(x,\xi)^{p}/p)$ for almost all $x\in\Omega$ and all $\xi\in\mathbb{R}^{n}$, $H$ is a family of norms on $\mathbb{R}^{n}$ ($n\geq 2$) parameterized by points $x\in\Omega$, and $V$ belongs to a local Morrey space. In particular, we investigate asymptotic behaviors of positive solutions of the equation near $\zeta$ and asymptotic behaviors of their quotients. 
\end{abstract}
\medskip

\noindent  \emph{2020 Mathematics Subject Classification.} Primary 35B53; Secondary 35B09, 35J92, 35B40.\\[2mm]
\noindent\emph{Key words and phrases.} Asymptotic behavior, Finsler $p$-Laplace equation, Fuchsian-type singularity, Morrey space, positive solutions.
\section{Introduction}
    Let~$n\geq 2$ be an integer, let~$1<p<\infty$, and let $\Omega\subseteq\mathbb{R}^{n}$ be a domain (i.e., a nonempty connected open set).
    We begin with four classical results, which are called the \emph{positive Liouville theorem}, the \emph{Picard principle}, the \emph{Poisson principle}, and the \emph{Riemann removable singularity theorem}, respectively.
    \btheorem[{see, e.g., \cite[p.~314]{Pin94}}]
    \begin{itemize}
    \item The set of all positive harmonic functions in $\R^{n}$ is a one-dimensional cone.
    \item The cone of all positive harmonic functions in the punctured unit ball in $\R^{n}$ vanishing on the unit sphere is one-dimensional.
    \item The cone of all positive harmonic functions in the unit ball in $\R^{n}$ vanishing on $S^{n-1}\setminus\{\xi\}$ is one-dimensional, where $S^{n-1}$ is the unit sphere in $\R^{n}$ and the point $\xi\in S^{n-1}$.
    \item If $u$ is a bounded harmonic function in a punctured neighborhood of the origin in $\R^{n}$, then $\myd{\lim_{x\rightarrow 0}u(x)}$ exists.
    \end{itemize}
    \etheorem
    In \cite{Pin94}, Pinchover studied positive Liouville-type theorems for a linear elliptic
operator of Fuchsian type and the asymptotic behavior of quotients of positive solutions of linear and semilinear equations. See also references therein. In \cite{Fraas1}, Fraas and Pinchover investigated positive Liouville-type theorems and the asymptotic behavior of positive solutions of the $p$-Laplace equation with a potential term. In \cite{Giri1,Giri2}, Giri and Pinchover extended results in \cite{Fraas1} to the~$(p,A)$-Laplace equation with a potential in a local Morrey space and a Wolff class, respectively. 

\textbf{Throughout this paper, we assume that $\zeta\in\{0,\infty\}$ is an isolated point of the boundary $\widehat{\partial}\Omega$ of $\Omega$ in the one-point compactification $\widehat{\R^{n}}$ of $\R^{n}$ with the ideal point $\infty$}. For every~$R>0$, let $I_{\zeta,R}\triangleq\begin{cases}
(0,R)& \mbox{if}~\zeta=0,\\
(R,\infty)& \mbox{if}~\zeta=\infty,	
\end{cases}$ and let $\Gamma_{\zeta,R}\triangleq\begin{cases}
B_{R}(0)\setminus\{0\}& \mbox{if}~\zeta=0,\\
\R^{n}\setminus D_{R}(0)& \mbox{if}~\zeta=\infty,\end{cases}$ where~$B_{R}(0)$ and $D_{R}(0)$ are the open and closed balls centered at $0$ with radius $R$, respectively. Clearly, there exists~$\hat{R}>0$ such that~$\Gamma_{\zeta,R}\subseteq\Omega$ for all~$R\in I_{\zeta,\hat{R}}$. We may also easily see that~$\Gamma_{\zeta,\hat{R}}\subseteq\Omega$. \textbf{Throughout this paper, we always fix such a positive number~$\hat{R}$, and assume that $\{H(x,\cdot)\}_{x\in\Gw}$ is a family of norms on $\mathbb{R}^{n}$ satisfying Assumptions~\ref{ass9} and that~$V\in M^{q}_{\loc}(p;\Omega)$, where $M^{q}_{\loc}(p;\Omega)$ is a local Morrey space defined in \Cref{Morreydef1}}. For any $\lambda\in\R$, let $\lambda^{<p-1>}\triangleq\lambda|\lambda|^{p-2}$. Under some further conditions, we study Fuchsian-type singularity (see \Cref{Fuchsian}) at $\zeta$ with respect to a norm $N\in C^{1}(\R^{n}\setminus\{0\})$ on $\R^{n}$ for the \emph{Finsler $p$-Laplace equation with a potential} 
 $$Q'[u]\triangleq Q'_{V}[u]\triangleq-\mathrm{div}\mathcal{A}(x,\nabla u)+Vu^{<p-1>}=0\qquad \mbox{ in } \Omega,$$
where $\mathcal{A}(x,\xi)\triangleq\nabla_{\xi}(H(x,\xi)^{p}/p)$ for almost all $x\in\Omega$ and all $\xi\in\mathbb{R}^{n}$. The equation is understood in the weak sense (see \Cref{def_sol}). Since we are are considering asymptotic behaviors near $\zeta$, it is natural to study positive solutions in $\Gamma_{\zeta,\hat{R}}$. For~$V=0$, the equation is called the \emph{Finsler $p$-Laplace equation} \cite{Combete}. For more on the Finsler $p$-Laplace equation, see, e.g., \cite{Bianchini,Jaros,Xianew}. The equation $Q'[u]=0$ was studied in \cite{Hou2,Hou3,HPR}. For example, in \cite{HPR}, criticality theory for the equation was established. Our main aim is to study asymptotic behaviors near~$\zeta$ of positive solutions of the equation with a Fuchsian-type singularity at~$\zeta$ and asymptotic behaviors of their quotients. We also prove a positive Liouville-type theorem for $p\leq n$.

Following \cite[Definition 3.1]{Giri1}, we define the Fuchsian-type singularity at $\zeta$ with respect to a norm $N\in C^{1}(\R^{n}\setminus\{0\})$ on $\R^{n}$. For every~$0<r<R$, let $\mathbb{A}_{r,R,N}\triangleq\{x\in\R^{n}~|~r<N(x)<R\}$. 
\bdefinition[Fuchsian-type singularity]\label{Fuchsian}
\emph{Let~$N\in C^{1}(\R^{n}\setminus\{0\})$ be a norm on~$\R^{n}$. The operator $Q'$ has a \emph{Fuchsian-type singularity} at~$\zeta$ with respect to $N$ if
\begin{itemize}
    \item there exist two positive constants $\alpha_{\zeta,\hat{R}}$ and~$\beta_{\zeta,\hat{R}}$ such that for almost all~$x\in\Gamma_{\zeta,\hat{R}}$ and all~$\xi\in\R^{n}$,
   \begin{align}\label{normereq}
\alpha_{\zeta,\hat{R}}|\xi|^{p}\le\mathcal{A}(x,\xi)\cdot\xi\quad\mbox{and}\quad
						|\mathcal{A}(x,\xi)|\le \beta_{\zeta,\hat{R}}\,|\xi|^{{p}-1};
   \end{align}
    \item there exist two positive numbers $a_{N}<1<b_{N}$ and a positive constant~$R_{N}$ such that the following function is bounded on~$I_{\zeta,R_{N}}$:\begin{align*}
\Lambda(R)\triangleq\begin{cases}
\|N(x)^{p-n/q}V\|_{M^{q}(p;\mathbb{A}_{a_{N}R,b_{N}R,N})}&\mbox{if}~p\neq n,\\
\|V\|_{M^{q}(n;\mathbb{A}_{a_{N}R,b_{N}R,N})}&\mbox{if}~p=n,
\end{cases}
\end{align*}
where~$\mathbb{A}_{a_{N}R,b_{N}R,N}\Subset\Gamma_{\zeta,\hat{R}}$ for all~$R\in I_{\zeta,R_{N}}$.
\end{itemize}}
\edefinition
\bexample[{\cite[Example 3.4]{Giri1}}]\label{pfuchsianexample}
\emph{Let~$\Omega=\R^{n}\setminus\{0\}$, and for almost all~$x\in\R^{n}\setminus\{0\}$, let~$H(x,\cdot)=|\cdot|$ be the Euclidean norm and let~$V(x)=\lambda/|x|^{p}$ ($\lambda\in\R$). Then the operator
$$Q'[u]=-\dive(|\nabla u|^{p-2}\nabla u)+\frac{\lambda u^{<p-1>}}{|x|^{p}},$$ has a Fuchsian-type singularity at both~$0$ and~$\infty$ with respect to $|\cdot|$, where $a_{|\cdot|}=1/2$ and $b_{|\cdot|}=3/2$.}
\eexample
\bexample\label{illus}
\emph{Let~$\Omega=\R^{n}\setminus\{0\}$ and let~$V(x)=\lambda/|x|^{p}$ ($\lambda\in\R$) for almost all~$x\in\R^{n}\setminus\{0\}$. If \eqref{normereq} holds, then by \Cref{pfuchsianexample}, the operator
$$Q'[u]=-\dive\mathcal{A}(x,\nabla u)+\frac{\lambda u^{<p-1>}}{|x|^{p}},$$ has a Fuchsian-type singularity at both~$0$ and~$\infty$ with respect to $|\cdot|$. For $p\neq n$ and $0<-\lambda\leq |(p-n)/p|^{p}$, $\lambda/|x|^{p}$ is called a \emph{Hardy potential}.}
\eexample
\bdefinition
\emph{The operator $Q'$ is called \emph{nonnegative} in $\Omega$ if the equation $Q'[u]=0$ has a positive (super-)solution in $\Omega$.}
\edefinition
\bdefinition
\emph{Suppose that~$Q'$ is nonnegative in~$\Omega$. Let~$\mathcal{G}_{\zeta}$ be the set of all positive solutions of $Q'[u]=0$ in connected punctured neighborhoods~$\Omega'\subseteq\Omega$ of~$\zeta$. The point~$\zeta$ is regular with respect to~$Q'$ if for all~$u,v\in\mathcal{G}_{\zeta}$,$$\lim_{\substack{x\rightarrow\zeta}}\frac{u(x)}{v(x)}\quad\mbox{exists in the wide sense.}$$}
\edefinition
\bdefinition
\emph{Suppose that~$Q'$ is nonegative in~$\Gamma_{\zeta,\hat{R}}$. Let~$\mathcal{G}_{\zeta,\hat{R}}$ be the set of all positive solutions of the equation~$Q'[u]=0$ in~$\Gamma_{\zeta,\hat{R}}$.}
\edefinition
When $H$ is independent of $x\in\Omega$, we say that $H$ is a fixed norm.
\bassumption\label{C1alphan}
\emph{Let $H$ be a fixed norm. Assume that~$\mathcal{A}(\cdot)\in C^{1}(\R^{n}\setminus\{0\})$ and that there exist two positive constants~$\Lambda_{1},\Lambda_{2}$ such that for all~$\xi\in\R^{n}\setminus\{0\}$ and all~$\eta\in\R^{n}$,
\begin{align*}
D\mathcal{A}(\xi)\eta\cdot\eta\geq \Lambda_{1}|\xi|^{p-2}|\eta|^{2}\quad\mbox{and}\quad|D\mathcal{A}(\xi)|\leq \Lambda_{2}|\xi|^{p-2},
\end{align*}
where the matrix function $D\mathcal{A}(\xi)$ is the derivative of $\mathcal{A}$.
}
\eassumption
We present one of our main results which can be easily concluded from \Cref{asympn}. 
\bcorollary\label{pncoro}
Let~$p\leq n$, let~$\zeta=0$, and let $H$ be a fixed norm satisfying \Cref{C1alphan}. Then~$0$ is regular with respect to~$Q'_{0}$.
\ecorollary

The paper is organized as follows. In Section 2, we define the Finsler $p$-Laplace equation with a potential and review some definitions and a result of criticality theory. In Section 3, we prove a Harnack convergence principle. In Section 4, we study a special class $\mathbf{M}_{\zeta,\hat{R}}$ of positive solutions of minimal growth near $S_{\hat{R}}(0)$, where $S_{\hat{R}}(0)$ is the sphere centered at $0$ with radius $\hat{R}$, and consider asymptotic behaviors of quotients of positive solutions of $Q'[w]=0$ in $\Gamma_{\zeta,\hat{R}}$. In Section 5, for a fixed norm $H$, we discuss asymptotic behaviors of positive Finsler $p$-harmonic functions and their quotients. In Section 6, we introduce a limiting dilated equation (operator). For $p\leq n$, we establish a positive Liouville-type theorem. 

\section*{Basic symbols}
 \normalsize
 \vspace{-1mm}
\begin{longtable}[1]{p{60pt} p{350pt} }
 $\mathbb{N}$ & the set of positive integers\\
$s'$ &the conjugate exponent~$s/(s-1)$ of~$s\in(1,\infty)$\\
$\zeta$& either~$0$ or~$\infty$ in~$\widehat{\R^{n}}$, $|\zeta|\triangleq\begin{cases}
		0&\mbox{if}~\zeta=0,\\	\infty&\mbox{if}~\zeta=\infty
		\end{cases}$\\
 $E/R$& $\{x/R\in\R^{n}~|~x\in E\}$ ($R>0$ and $E\subseteq\R^{n}$)\\
 $\mathbb{A}_{r,R}$& $\{x\in\R^{n}~|~r<|x|<R\}$ ($r,R>0$)\\
 $B_r(x)$  & the open ball centered at $x\in\mathbb{R}^{n}$ with radius $r>0$\\
$S_r(x)$  & the sphere centered at $x\in\mathbb{R}^{n}$ with radius $r>0$\\
$D_r(x)$  & the closed ball centered at $x\in\mathbb{R}^{n}$ with radius $r>0$\\
$B_{r,N}$  & the open ball with respect to a norm~$N$ on~$\R^{n}$ centered at $0$ with radius $r>0$, i.e., $\{x\in\R^{n}~|~N(x)<r\}$\\
$S_{r,N}$  & the sphere with respect to a norm~$N$ on~$\R^{n}$ centered at $0$ with radius $r>0$, i.e., $\{x\in\R^{n}~|~N(x)=r\}$\\
$D_{r,N}$  & the closed ball with respect to a norm~$N$ on~$\R^{n}$ centered at $0$ with radius $r>0$, i.e., $\{x\in\R^{n}~|~N(x)\leq r\}$\\
  $C(\Omega)$& the space of continuous functions in~$\Omega$\\
  $E \Subset F$& $\overline{E}$ is a compact subset of $F$ with~$E,F\subseteq\R^{n}$\\
 $\vol(E)$ & the Lebesgue measure of a measurable subset~$E$ of~$\R^{n}$\\
     $\sgn(\cdot)$& the sign function on~$\R$\\ 
 $\mathring{K},\overline{K}$& the interior and closure of~$K\subseteq\R^{n}$, respectively \\
 $\mathrm{diam}(K)$ & the diameter of~$K\subseteq\R^{n}$\\
$\supp f$& the support of a function~$f$\\
$\mathcal{F}_{c}(U)$ & the space~$\{f\in\mathcal{F}(U)~|~\supp f\Subset U\}$, where~$U\subseteq\R^{n}$ is a nonempty open set and~$\mathcal{F}(U)$ is an arbitrary linear space of functions on~$U$\\
$C$ & a positive constant which may vary from line to line\\
$f\asymp g$&$cg\leq f\leq Cg$ for some positive constants~$c$ and~$C$, where $f$ and $g$ are nonnegative functions with the same domain\\
		\end{longtable}
\section{Finsler $p$-Laplace equation with a potential}
    In this section, we discuss the norm family~$H$, a variational Lagrangian~$F$, and the operator~$\mathcal{A}$, define the local Morrey space $M^{q}_{\loc}(p;\Omega)$, the energy functional $Q$, and the Finsler $p$-Laplace equation with a potential, and collect some definitions and a result in criticality theory. For more on criticality theory, see, e.g., \cite{Fischer1,Fischer2,Giri2,HPR,Keller2,Kovarik,PPAPDE,Regev,PCVPDE,VRobin2} and references therein.
   
   We start with the norm family~$H$. 
\begin{assumptions}\label{ass9} 
		{\em			Let~$H:\Omega\times\R^{n}$ be a mapping such that for almost all~$x\in\Omega$,~$H(x,\cdot)$ is a norm on $\mathbb{R}^{n}$, with the following properties:   
			\begin{itemize}
				\item {\bf (Measurability)} For all~$\xi\in\mathbb{R}^{n}$, the mapping $x\mapsto H(x,\xi)$ is measurable in $\Gw$;
				\item {\bf (Local uniform equivalence)} for every domain $\omega\Subset \Gw$, there exist two constants~$0<\kappa_\omega\leq\nu_\omega<\infty$ such that for almost all $x\in \omega$ and all $\xi\in \R^n$,$$\kappa_\omega |\xi|\leq H(x,\xi)\leq\nu_\omega|\xi|;$$
				\item {\bf (Uniform convexity)} for almost all~$x\in\Omega$, the Banach space $\mathbb{R}^{n}_{x}\triangleq(\mathbb{R}^{n},H(x,\cdot))$ is uniformly convex; 
				\item {\bf (Differentiability with respect to $\xi$)} for almost all~$x\in \Gw$, the mapping $\xi\mapsto H(x,\xi)$ is differentiable in $\R^n\setminus\{0\}$.
			\end{itemize}		
		}
	\end{assumptions}
    The norm family $H$ satisfying \Cref{ass9} produces the following variational Lagrangian $F$.
    \begin{theorem}[{\cite[Theorem~2.4]{Hou2}}]
  For almost all $x\in \Omega$ and all~$\xi\in\R^{n}$, let
	$$F(x,\xi)\triangleq\frac{1}{p}H(x,\xi)^{p}.$$
            Then
			\begin{itemize}
				\item {\bf (Measurability)} for all~$\xi\in\mathbb{R}^{n}$, the mapping $x\mapsto F(x,\xi)$ is measurable in $\Gw$;
				\item {\bf (Local uniform ellipticity and boundedness)} 
				for every domain $\omega\Subset \Gw$, there exist two constants $0<\bar{\kappa}_\omega\leq\bar{\nu}_\omega<\infty$ such that for almost all $x\in \omega$ and all $\xi\in \R^n$,$$\bar{\kappa}_\omega|\xi|^{p}\leq F(x,\xi)\leq\bar{\nu}_\omega|\xi|^{p};$$
				\item {\bf (Strict convexity and $C^1$ with respect to $\xi$)}  for almost all~$x\in \Gw$, the mapping $\xi\mapsto F(x,\xi)$ is strictly convex and continuously differentiable in $\R^n$;
				\item {\bf (Homogeneity)} for almost all~$x\in \Gw$, all~$\lambda\in\mathbb{R}$, and all~$\xi\in\mathbb{R}^{n}$, $F(x,\lambda\xi)=|\lambda|^{p}F(x,\xi)$.
			\end{itemize}		
	\end{theorem}
    \begin{notation}
\emph{Recall that for almost all~$x\in\Omega$ and all~$\xi\in\mathbb{R}^{n}$, $$\mathcal{A}(x,\xi)=\nabla_{\xi}(H(x,\xi)^{p}/p)=\nabla_\xi F(x,\xi).$$ For almost all~$x\in\Omega$ and all~$\xi\in\mathbb{R}^{n}$, we also denote~$H(x,\xi)$ by~$|\xi|_{\mathcal{A}}$ in some places. For a fixed norm $H$, we simply write~$H(\xi)$ and~$\mathcal{A}(\xi)$.}
	\end{notation}
The homogeneity of~$F$ yields the following relationship between~$F$ and~$\mathcal{A}$.
    \begin{lemma}[{\cite[p.~100]{HKM}}]	For almost all~$x\in\Omega$ and all~$\xi\in\mathbb{R}^{n}$,~$\mathcal{A}(x,\xi)\cdot\xi =pF(x,\xi).$
	\end{lemma}
Some basic properties of the operator~$\mathcal{A}$ are listed in the following theorem.
    \begin{Thm}[{\cite[Lemma~5.9]{HKM}}]\label{thm_1}
			For every domain $\omega\Subset\Omega$, let  $\alpha_{\omega}=\bar{\kappa}_{\omega}$, $\beta_{\omega}=2^{p}\bar{\nu}_{\omega}$. Then the vector-valued function~$\mathcal{A}:  \Gw\times \mathbb{R}^{n}\rightarrow \mathbb{R}^{n}$ satisfies the following properties:
				\begin{itemize}
					\item {\bf (Continuity and measurability)} For almost all $x\in \Gw$, the function $\xi\mapsto \mathcal{A}(x,\xi )$ is continuous in~$\R^{n}$, and for all~$\xi\in \mathbb{R}^{n}$, the function $x \mapsto \mathcal{A}(x,\xi)$ is Lebesgue measurable on $\Gw$;
					
					%
					\item {\bf (Local uniform ellipticity and boundedness)} for all domains $\omega\Subset \Gw$, almost all $x\in \omega$, and all $\xi \in \mathbb{R}^{n}$,		\begin{equation*}\label{structure}
						\alpha_\omega|\xi|^{p}\le\mathcal{A}(x,\xi)\cdot\xi\quad\mbox{and}\quad
						|\mathcal{A}(x,\xi)|\le \beta_\omega\,|\xi|^{{p}-1};
					\end{equation*}
					\item {\bf (Strict monotonicity)} for almost all~$x\in\Gw$ and all~$\xi,\eta\in\mathbb{R}^{n}$ ($\xi\neq\eta$),
					$$\big(\mathcal{A}(x,\xi)-\mathcal{A}(x,\eta)\big) \! \cdot \! (\xi-\eta)> 0;$$
                    \item {\bf (Homogeneity)} for all~$\lambda\in {\mathbb{R}\setminus\{0\}}$,
					$\mathcal{A}(x,\lambda \xi)=\lambda^{<p-1>}\mathcal{A}(x,\xi)$.
			\end{itemize}
		\end{Thm}
         \begin{Def}\label{Morreydef1}
  \emph{Suppose that~$\omega\Subset\Omega$ is a domain and that $f$ is a real-valued measurable function on~$\omega$. 
				\begin{itemize}
					\item For $p<n$ and $q>n/p$, we define$$\Vert f\Vert_{M^{q}(p;\,\omega)}\triangleq \sup_{\substack{y\in\gw\\0<r<\diam(\gw)}}
					\frac{1}{r^{n/q'}}\int_{\omega\cap B_{r}(y)}|f|\dx$$ and $$M^{q}(p;\omega)\triangleq\{f\in L^{1}_{\loc}(\omega)~|~\Vert f\Vert_{M^{q}(p;\,\omega)}<\infty\};$$ 
					\item for $p=n$ and $q>n$ , we define$$\Vert f\Vert_{M^{q}(n;\,\omega)}\triangleq \sup_{\substack{y\in\gw\\0<r<\diam(\gw)}} \varphi_{q}(r)\int_{\omega\cap B_{r}(y)}|f|\dx,$$
					where $\varphi_{q}(r)\triangleq \left(\log\big(\mathrm{diam}(\omega)/r\big)\right)^{q/n'}$, and
$$M^{q}(n;\omega)\triangleq\{f\in L^{1}_{\loc}(\omega)~|~\Vert f\Vert_{M^{q}(n;\,\omega)}<\infty\};$$
					\item for $p>n$ and $q=1$, we define~$M^{q}(p;\omega)\triangleq L^{1}(\omega)$.
				\end{itemize}
    Finally, we define the \emph{local Morrey space} by $$M^{q}_{\loc}(p;\Omega)\triangleq \bigcap_{\substack{\omega\Subset\Omega\\\omega~\mbox{is a domain}}}M^{q}(p;\omega).$$
			}
		\end{Def}
  \bdefinition
\emph{Let $\omega\Subset\Omega$ be a domain. For $p<n$, the space~$\mathbb{M}^{q}(p;\omega)$ consists of all functions $f$ in~$M^{q}(p;\omega)$ such that~$|f|^{\left(\frac{p^{*}-1}{p-1}\right)'}\in M^{q}(p;\omega)$; for $p=n$, $\mathbb{M}^{q}(p;\omega)$ is defined as $L^{\rho}(\omega)$ for some~$\rho>1$; for $p>n$, $\mathbb{M}^{q}(p;\omega)$ is defined as $M^{q}(p;\omega)=L^{1}(\omega)$. Furthermore, we define $$\mathbb{M}^{q}_{\loc}(p;\Omega)\triangleq \bigcap_{\substack{\omega\Subset\Omega\\\omega~\mbox{is a domain}}}\mathbb{M}^{q}(p;\omega).$$}
                    \edefinition
  \begin{definition}\label{efdfn}
 \emph{The \emph{energy functional}~$Q$ on~$\core$ is defined by \begin{align*}
Q[\phi]\triangleq Q[\phi;\Omega]\triangleq Q_{\mathcal{A},V}[\phi]\triangleq Q_{V}[\phi]\triangleq\int_{\Omega}\left(H(x,\nabla \phi)^{p}+ V|\phi|^{p}\right)\dx.
\end{align*}The functional~$Q$ is called \emph{nonnegative} in~$\Omega$ if $Q[\phi]\geq 0$ for all~$\phi\in C^{\infty}_{c}(\Omega)$.}
 \end{definition} 
  
 Now we define the \emph{nonhomogeneous Finsler $p$-Laplace equation with a potential}.
 \begin{Def}\label{def_sol}
			\emph{  
Let~$g\in M^{q}_{\loc}(p;\Omega)$. A function~$v\in W^{1,p}_{\loc}(\Omega)$ is a  (\emph{weak}) \emph{solution} of the equation
				\begin{equation*}\label{half}
					Q'[u]\triangleq Q'_{\mathcal{A},V}[u]\triangleq Q'_{V}[u]\triangleq -\dive\mathcal{A}(x,\nabla u)+Vu^{<p-1>}=g,
				\end{equation*}
				in~$\Omega$ if for all~$\phi \in C_{c}^{\infty}(\Omega)$,$$\int_{\Omega}\mathcal{A}(x,\nabla v)\cdot \nabla \phi\dx+\int_{\Omega}Vv^{<p-1>} \phi\dx=\int_{\Omega}g\phi\dx.$$
    A function~$v\in W^{1,p}_{\loc}(\Omega)$ is a (\emph{weak}) \emph{supersolution (subsolution)} of the above equation
				in~$\Omega$ if for all nonnegative~$\phi \in C_{c}^{\infty}(\Omega)$, $$\int_{\Omega}\mathcal{A}(x,\nabla v)\cdot \nabla \phi\dx+\int_{\Omega}Vv^{<p-1>}\phi\dx\geq (\leq) \int_{\Omega}g\phi\dx.$$ 
                A solution of~$Q'_{0}[u]=0$ in $\Omega$ is called \emph{Finsler $p$-harmonic} in $\Omega$.}
                \end{Def}
\begin{Def}\label{gpeigen}
	{\em
		The \emph{principal eigenvalue} of $Q'$ in a domain $\omega\Subset \Gw$
		is defined by
		$$\lambda_{1}(Q;\omega)\triangleq\inf_{\phi\in C^{\infty}_{c}(\omega) \setminus\{0\}}\frac{Q[\phi;\omega]}{\Vert \phi\Vert_{L^{p}(\omega)}^{p}}.$$}	
\end{Def}
\bassumption\label{C1alpha}
\emph{For all~$x\in\Omega$,~$\mathcal{A}(x,\cdot)\in C^{1}(\R^{n}\setminus\{0\})$. For every domain~$\omega\Subset\Omega$, 
there exist positive constants~$\vartheta\leq 1, \Lambda_{1},\Lambda_{2}$, and~$\Lambda_{3}$ such that for all~$x,y\in\omega$, all~$\xi\in\R^{n}\setminus\{0\}$, and all~$\eta\in\R^{n}$,
\begin{align*}
D_{\xi}\mathcal{A}(x,\xi)\eta\cdot\eta\geq \Lambda_{1}|\xi|^{p-2}|\eta|^{2},\quad |D_{\xi}\mathcal{A}(x,\xi)|\leq \Lambda_{2}|\xi|^{p-2},
\end{align*}
and \begin{align*}
|\mathcal{A}(x,\xi)-\mathcal{A}(y,\xi)|\leq \Lambda_{3}|\xi|^{p-1}|x-y|^{\vartheta}.
\end{align*}}
\eassumption
\begin{lemma}[{\cite[Lemma 2.10]{Hou3}}]\label{bllem}
	{\em
 Suppose that \Cref{C1alpha} holds. For every domain $\gw\Subset \Gw$, there exists  
		a positive constant $C$ (independent of~$x\in\omega$ and~$\xi,\eta\in\R^{n}$) such that for all~$x\in \gw$ and~$\xi,\eta\in \R^n$,  
	\begin{align*}	|\xi+\eta|^{p}_{\mathcal{A}}-|\xi|^{p}_{\mathcal{A}}-p\mathcal{A}(x,\xi)\cdot\eta\geq \begin{cases}
				C|\eta|^{p}&\mbox{if $p\geq 2$,}\\
				C|\eta|^{2}(|\xi|+|\eta|)^{p-2}&\mbox{if $p<2$.}
		\end{cases}
   \end{align*}	
	}
\end{lemma}
\bdefinition
\emph{The operator $Q'$ is \emph{critical} in $\Omega$ if $Q'$ is nonnegative in $\Omega$ and $Q'_{V-W}$ is not nonnegative in $\Omega$ for any nonnegative $W\in M^{q}_{\loc}(p;\Omega)\setminus\{0\}$.}
\edefinition
By \cite[Theorems 5.3, 6.9, and 6.12]{HPR}, we may get the following result. 
\blemma\label{uniques}
Suppose that \Cref{C1alpha} holds in $\Omega$ and that $q>n$ for $p<2$. Suppose that $Q'$ is critical in $\Omega$. Then $Q'[u]=0$ has a unique (up to a positive multiplicative constant) positive supersolution (in $W^{1,\infty}_{\loc}(\Omega)$ for $p<2$) in $\Omega$, which is also a unique (up to a positive multiplicative constant) positive solution of the equation $Q'[u]=0$ in $\Omega$.
\elemma
\bremark
\emph{Up to a positive multiplicative constant, the unique positive supersolution (in $W^{1,\infty}_{\loc}(\Omega)$ for $p<2$) in $\Omega$, which is actually a solution, is called an \emph{Agmon ground state}}.
\eremark
\begin{definition}
\emph{A compact subset~$K$ of~$\R^{n}$ is called \emph{admissible} if~$K$ is the closure of a Lipschitz domain with connected boundary.}
\end{definition}
\begin{Def}\label{dfnmg}
\emph{Let $K_{0}$ be a compact subset of~$\Omega$ such that~$\Omega\setminus K_{0}$ is a domain. A positive solution~$u$ of~$Q'[w]=0$ in~$\Omega\setminus K_{0}$ has \emph{minimal growth in a neighborhood of infinity} in $\Omega$ if for all admissible compact subsets~$K$ of~$\Omega$ with~$K_{0}\subseteq \mathring{K}$, and all positive solutions~$v\in C\left(\Omega\setminus \mathring{K}\right)$ of~$Q'[w]=g$  in~$\Omega\setminus K$ such that $u\leq v$ on~$\partial K$, it holds that~$ u\leq v$ in~$\Omega\setminus K$, where~$g\in M^{q}_{\loc}(p;\Omega)$ for $p\neq n$, $g\in L^{\bar{\rho}}_{\loc}(\Omega)$ ($\bar{\rho}>1$) for $p=n$, and $g$ is nonnegative a.e. in $\Omega\setminus K$. The set of all such positive solutions is denoted by~$\mathcal{M}_{\Omega;K_{0}}$. If~$K_{0}=\emptyset$, then $u$ is called a \emph{global minimal positive solution} of~$Q'[w]=0$ in~$\Omega$.}
  \end{Def}
\section{Harnack convergence principle}
In this section, we prove a Harnack convergence principle.
\bdefinition
\emph{A \emph{smooth (Lipschitz) exhaustion} of~$\Omega$ is a sequence of~$C^{\infty}$ (Lipschitz) domains~$\{\omega_{k}\}_{k\in\mathbb{N}}$ such that for all~$k\in\mathbb{N}$,~$\omega_{k}\Subset\omega_{k+1}\Subset\Omega$ and~$\cup_{k\in\mathbb{N}}\omega_{k}=\Omega$.}
\edefinition
\btheorem [Harnack convergence principle]\label{hcpg}
Let $\{\omega_{k}\}_{k\in\mathbb{N}}$ be a Lipschitz exhaustion of~$\Omega$. Let $\{\mathcal{A}_{k}(x,\cdot)\}_{x\in\omega_{k}}$ be a sequence of operators with all the properties in \Cref{thm_1} in~$\omega_{k}$ for every~$k\in\mathbb{N}$ such that for every domain compactly contained in $\Omega$ and all sufficiently large~$k\in\mathbb{N}$, the two constants in the local uniform ellipticity and boundedness conditions are independent of~$k\in\mathbb{N}$, that~$\mathcal{A}_{k}(x,\cdot)$ is pointwise equicontinuous in~$\R^{n}$ for almost every~$x\in\Omega$ and all sufficiently large~$k\in\mathbb{N}$, and that~$\lim_{k\rightarrow\infty}\mathcal{A}_{k}(x,\xi)=\mathcal{A}(x,\xi)$ for almost all~$x\in\Omega$ and all~$\xi\in\R^{n}$. Suppose that $V_{k}\in M^{q}_{\loc}(p;\omega_{k})$ converges to~$V\in M^{q}_{\loc}(p;\Omega)$ weakly in~$M^{q}_{\loc}(p;\Omega)$ as~$k\rightarrow\infty$. For every~$k\in\mathbb{N}$, let $v_{k}$ be a positive solution of $Q'_{\mathcal{A}_{k},V_{k}}[u]=0$ in~$\omega_{k}$ such that for all~$k\in\mathbb{N}$ and some~$x_{0}\in\omega_{1}$,~$v_{k}(x_{0})=1$. Then up to a subsequence,~$\{v_{k}\}_{k\in\mathbb{N}}$ converges weakly in~$W^{1,p}_{\loc}(\Omega)$ and locally uniformly in~$\Omega$ to a positive solution of~$Q'[u]=0$ in~$\Omega$.
\etheorem
\bproof
The proof is similar to that of \cite[Proposition 2.7]{Giri1}.

We claim that for every domain~$\omega\Subset\Omega$ such that~$x_{0}\in\omega$,~$\{v_{k}\}_{k\in\mathbb{N}}$ has a subsequence converging uniformly in~$\omega$ to a H\"older continuous function~$v_{\omega}$. 

By the local Schauder-type estimate \cite[Theorem 3.2]{HPR}, for every~$k\in\mathbb{N}$,~$v_{k}$ is locally H\"{o}lder continuous in~$\omega_{k}$. Since~$V_{k}\in M^{q}_{\loc}(p;\omega_{k})$ converges to~$V\in M^{q}_{\loc}(p;\Omega)$ weakly in~$M^{q}_{\loc}(p;\Omega)$ as~$k\rightarrow\infty$, for every domain~$\omega\Subset\Omega$,~$\Vert V_{k}\Vert_{M^{q}(p;\,\omega)}$ is bounded for all sufficiently large~$k\in\mathbb{N}$ (see, e.g., \cite[p.~66 and Theorem~3.18]{Rudin}). Because~$v_{k}(x_{0})=1$ for all~$k\in\mathbb{N}$ and the two constants in the local uniform ellipticity and boundedness conditions are independent of~$k\in\mathbb{N}$ for every domain $\omega\Subset\Omega$ and all sufficiently large~$k\in\mathbb{N}$, by the local Harnack inequality \cite[Theorem 3.1]{HPR}, for every domain~$\omega\Subset\Omega$ such that~$x_{0}\in\omega$,~$v_{k}$ is uniformly bounded in~$\omega$ with respect to all sufficiently large~$k\in\mathbb{N}$. By the local Schauder-type estimate \cite[Theorem 3.2]{HPR}, for every domain~$\omega\Subset\Omega$ such that~$x_{0}\in\omega$,~$v_{k}$ is equicontinuous in~$\omega$ for all sufficiently large~$k\in\mathbb{N}$. By the Arzel\`a--Ascoli theorem \cite[Theorem 11.28]{RudinRC}, we get the claim.

Let~$\omega\Subset\Omega$ be an arbitrary domain and we fix a positive integer~$k_{0}$ such that~$\omega\Subset\omega_{k_{0}}$. We claim that with respect to all sufficiently large positive integers~$k$,~$v_{k}$ is bounded in~$W^{1,p}(\omega)$ and that up to a subsequence,~$\{v_{k}\}_{k\in\mathbb{N}}$ converges weakly in~$W^{1,p}_{\loc}(\Omega)$ and locally uniformly in~$\Omega$ to some nonnegative function~$v\in W^{1,p}_{\loc}(\Omega)\cap C(\Omega)$, where~$v(x_{0})=1$.

For all~$\phi\in C^{\infty}_{c}(\omega_{k_{0}})$ and all positive integers~$k>k_{0}$,~$v_{k}|\phi|^{p}\in W^{1,p}_{c}(\omega_{k_{0}})$. Since~$v_{k}$ is a positive solution of~$Q'_{\mathcal{A}_{k},V_{k}}[u]=0$ in~$\omega_{k_{0}}$, we have
\begin{align*}
0&=\int_{\omega_{k_{0}}}\mathcal{A}_{k}(x,\nabla v_{k})\cdot\nabla(v_{k}|\phi|^{p})dx+\int_{\omega_{k_{0}}}V_{k}v_{k}^{p-1}v_{k}|\phi|^{p}dx\\
&=\int_{\omega_{k_{0}}}\mathcal{A}_{k}(x,\nabla v_{k})\cdot\left(|\phi|^{p}\nabla v_{k}+pv_{k}|\phi|^{p-1}\nabla|\phi|\right)dx+\int_{\omega_{k_{0}}}V_{k}v_{k}^{p}|\phi|^{p}dx.
\end{align*}
Next we denote the two constants in the local uniform ellipticity and boundedness of~$\mathcal{A}_{k}$ in~$\omega_{k_{0}}$ by~$\bar{\alpha}_{\omega_{k_{0}}}$ and $\bar{\beta}_{\omega_{k_{0}}}$ for all sufficiently large~$k\in\mathbb{N}$, respectively. Then for all sufficiently large~$k\in\mathbb{N}$,
\begin{align*}
&\int_{\omega_{k_{0}}}|\nabla v_{k}|_{\mathcal{A}_{k}}^{p}|\phi|^{p}dx\\
&\leq p\int_{\omega_{k_{0}}}|\phi|^{p-1}v_{k}|\mathcal{A}_{k}(x,\nabla v_{k})\cdot\nabla|\phi||dx+\int_{\omega_{k_{0}}}|V_{k}|v_{k}^{p}|\phi|^{p}dx\\
&\leq p\bar{\beta}_{\omega_{k_{0}}}\int_{\omega_{k_{0}}}|\phi|^{p-1}v_{k}|\nabla v_{k}|^{p-1}|\nabla|\phi||dx+\int_{\omega_{k_{0}}}|V_{k}|v_{k}^{p}|\phi|^{p}dx\\
&=p\bar{\beta}_{\omega_{k_{0}}}\int_{\omega_{k_{0}}}v_{k}|\nabla\phi||\phi|^{p-1}|\nabla v_{k}|^{p-1}dx+\int_{\omega_{k_{0}}}|V_{k}|v_{k}^{p}|\phi|^{p}dx\quad(|\nabla|\phi||=|\nabla\phi|~\mbox{in}~\omega_{k_{0}})\\
&\leq \varepsilon\bar{\beta}_{\omega_{k_{0}}}^{p'}\int_{\omega_{k_{0}}}|\phi|^{p}|\nabla v_{k}|^{p}dx+\left(\frac{p-1}{\varepsilon}\right)^{p-1}\int_{\omega_{k_{0}}}v_{k}^{p}|\nabla\phi|^{p}dx\\
&+\delta\int_{\omega_{k_{0}}}|\nabla(v_{k}\phi)|^{p}dx+C\int_{\omega_{k_{0}}}|v_{k}\phi|^{p}dx\\
&(\mbox{by Young's inequality and the Morrey-Adams theorem \cite[Theorem 2.4 (2)]{PPAPDE}}),
\end{align*}
where~$\varepsilon>0$,~$0<\delta\leq\delta_{0}$, and~$C=C\left(n,p,q,\omega_{k_{0}},\omega_{k_{0}+1},\delta,\Vert V_{k}\Vert_{M^{q}(p;\,\omega_{k_{0}+1})}\right)$.  
Furthermore, for all sufficiently large positive integers~$k$, all~$0<\varepsilon<\bar{\alpha}_{\omega_{k_{0}}}/\bar{\beta}_{\omega_{k_{0}}}^{p'}$, all~$0<\delta\leq\delta_{0}$, and~$C=C\left(n,p,q,\omega_{k_{0}},\omega_{k_{0}+1},\delta,\Vert V_{k}\Vert_{M^{q}(p;\,\omega_{k_{0}+1})}\right)$, 
\begin{align*}
&(\bar{\alpha}_{\omega_{k_{0}}}-\varepsilon\bar{\beta}_{\omega_{k_{0}}}^{p'})\int_{\omega_{k_{0}}}|\nabla v_{k}|^{p}|\phi|^{p}dx\leq\int_{\omega_{k_{0}}}|\nabla v_{k}|_{\mathcal{A}_{k}}^{p}|\phi|^{p}dx-\varepsilon\bar{\beta}_{\omega_{k_{0}}}^{p'}\int_{\omega_{k_{0}}}|\nabla v_{k}|^{p}|\phi|^{p}dx\\
&\leq\left(\frac{p-1}{\varepsilon}\right)^{p-1}\int_{\omega_{k_{0}}}v_{k}^{p}|\nabla\phi|^{p}dx+\delta\int_{\omega_{k_{0}}}|\nabla(v_{k}\phi)|^{p}dx+C\int_{\omega_{k_{0}}}|v_{k}\phi|^{p}dx\\
&\leq\left(\frac{p-1}{\varepsilon}\right)^{p-1}\int_{\omega_{k_{0}}}v_{k}^{p}|\nabla\phi|^{p}dx+2^{p-1}\delta\int_{\omega_{k_{0}}}|v_{k}\nabla\phi|^{p}dx+2^{p-1}\delta\int_{\omega_{k_{0}}}|\phi\nabla v_{k}|^{p}dx\\
&+C\int_{\omega_{k_{0}}}|v_{k}\phi|^{p}dx.
\end{align*}
Therefore, for all sufficiently large positive integers~$k$, all~$0<\varepsilon<\bar{\alpha}_{\omega_{k_{0}}}/\bar{\beta}_{\omega_{k_{0}}}^{p'}$, all~$0<\delta\leq\delta_{0}$ such that~$\bar{\alpha}_{\omega_{k_{0}}}-\varepsilon\bar{\beta}_{\omega_{k_{0}}}^{p'}-2^{p-1}\delta>0$, and~$C=C\left(n,p,q,\omega_{k_{0}},\omega_{k_{0}+1},\delta,\Vert V_{k}\Vert_{M^{q}(p;\,\omega_{k_{0}+1})}\right)$,
\begin{align*}
&\left(\bar{\alpha}_{\omega_{k_{0}}}-\varepsilon\bar{\beta}_{\omega_{k_{0}}}^{p'}-2^{p-1}\delta\right)\int_{\omega_{k_{0}}}|\nabla v_{k}|^{p}|\phi|^{p}dx\\
&\leq\left(\frac{p-1}{\varepsilon}\right)^{p-1}\int_{\omega_{k_{0}}}v_{k}^{p}|\nabla\phi|^{p}dx+2^{p-1}\delta\int_{\omega_{k_{0}}}|v_{k}\nabla\phi|^{p}dx+C\int_{\omega_{k_{0}}}|v_{k}\phi|^{p}dx\\
&\leq\left(\left(\frac{p-1}{\varepsilon}\right)^{p-1}+2^{p-1}\delta\right)\int_{\omega_{k_{0}}}v_{k}^{p}|\nabla\phi|^{p}dx+C\int_{\omega_{k_{0}}}|v_{k}\phi|^{p}dx.
\end{align*}
Pick a domain~$\omega'$ such that~$\omega\Subset\omega'\Subset\omega_{k_{0}}$ and a function~$0\leq\phi\in C^{\infty}_{c}(\omega_{k_{0}})\leq 1$ such that~$\supp\phi\subseteq\omega'$ and that~$\phi|_{\omega}\equiv 1$. We deduce that for all sufficiently large positive integers~$k$, all~$0<\varepsilon<\bar{\alpha}_{\omega_{k_{0}}}/\bar{\beta}_{\omega_{k_{0}}}^{p'}$, all~$0<\delta\leq\delta_{0}$ such that~$\bar{\alpha}_{\omega_{k_{0}}}-\varepsilon\bar{\beta}_{\omega_{k_{0}}}^{p'}-2^{p-1}\delta>0$, and~$C=C\left(n,p,q,\omega_{k_{0}},\omega_{k_{0}+1},\delta,\Vert V_{k}\Vert_{M^{q}(p;\,\omega_{k_{0}+1})}\right)$,
\bals
&\int_{\omega}|\nabla v_{k}|^{p}\dx+\int_{\omega}|v_{k}|^{p}\dx\leq \int_{\omega_{k_{0}}}|\nabla v_{k}|^{p}\phi^{p}\dx+\int_{\omega_{k_{0}}}|v_{k}|^{p}\phi^{p}\dx\\
&\leq C'\int_{\omega_{k_{0}}}v_{k}^{p}|\nabla\phi|^{p}dx+C''\int_{\omega_{k_{0}}}|v_{k}\phi|^{p}dx,
\eals
where
$$C'=\frac{\left(\frac{p-1}{\varepsilon}\right)^{p-1}+2^{p-1}\delta}{\bar{\alpha}_{\omega_{k_{0}}}-\varepsilon\bar{\beta}_{\omega_{k_{0}}}^{p'}-2^{p-1}\delta}\quad\mbox{and}\quad C''=\frac{C}{\bar{\alpha}_{\omega_{k_{0}}}-\varepsilon\bar{\beta}_{\omega_{k_{0}}}^{p'}-2^{p-1}\delta}+1.$$
Recall that with respect to all sufficiently large positive integers~$k$,~$\Vert V_{k}\Vert_{M^{q}(p;\,\omega_{k_{0}+1})}$ is bounded and~$v_{k}$ is uniformly bounded in~$\omega_{k_{0}}$. The boundedness claim follows. Recall that up to a subsequence,~$\{v_{k}\}_{k\in\mathbb{N}}$ converges uniformly in~$\omega_{k_{0}}$ to a H\"older continuous function~$v_{\omega_{k_{0}}}$. Thus, up to a subsequence,~$\{v_{k}\}_{k\in\mathbb{N}}$ converges weakly in~$W^{1,p}(\omega)$ and uniformly in~$\omega$ to~$v_{\omega}\in W^{1,p}(\omega)\cap C^{\gamma}(\omega)$ ($\gamma\in(0,1])$. Finally, by a diagonal argument, we also get the convergence claim. 

We claim that the function~$v$ is a solution of~$Q'[u]=0$ in~$\Omega$. 

We prove this by showing that for every~$k\in\mathbb{N}$ and all~$\psi\in C^{\infty}_{c}(\omega_{k})$,$$\int_{\omega_{k}}\mathcal{A}(x,\nabla v)\cdot \nabla \psi\dx+\int_{\omega_{k}}Vv^{p-1}\psi\dx=0.$$ Next fix an arbitrary~$k\in\mathbb{N}$. Note that for all positive integers~$i\geq k$ and all~$\psi\in C^{\infty}_{c}(\omega_{k})$,
$$\int_{\omega_{k}}\mathcal{A}_{i}(x,\nabla v_{i})\cdot \nabla \psi\dx+\int_{\omega_{k}}V_{i}v_{i}^{p-1}\psi\dx=0.$$ For all~$\psi\in C^{\infty}_{c}(\omega_{k})$, the same proof as the counterpart of \cite[Theorem 3.5]{HPR} yields, up to a subsequence,
\bals
\lim_{i\rightarrow\infty}\int_{\omega_{k}}V_{i}v_{i}^{p-1}\psi\dx=\int_{\omega_{k}}Vv^{p-1}\psi\dx.
\eals
It suffices to prove that up to a subsequence, for all~$\psi\in C^{\infty}_{c}(\omega_{k})$,
\bals
\lim_{i\rightarrow\infty}\int_{\omega_{k}}\mathcal{A}_{i}(x,\nabla v_{i})\cdot \nabla \psi\dx=\int_{\omega_{k}}\mathcal{A}(x,\nabla v)\cdot \nabla \psi\dx.
\eals
To this end, we only need to show that~$\{\mathcal{A}_{i}(x,\nabla v_{i})\}_{i\in\mathbb{N}}$ converges to~$\mathcal{A}(x,\nabla v)$ weakly in~$L^{p'}(\omega_{k};\R^{n})$ up to a subsequence. Recall that for all sufficiently large~$i\in\mathbb{N}$, the two constants in the local uniform ellipticity and boundedness conditions are independent of~$i\in\mathbb{N}$ and~$v_{i}$ is bounded in~$W^{1,p}(\omega_{k})$. Then~$\mathcal{A}_{i}(x,\nabla v_{i})$ is bounded in~$L^{p'}(\omega_{k};\R^{n})$ with respect to all sufficiently large~$i\in\mathbb{N}$.
By \cite[Theorem~13.44]{Hewitt}, the task boils down to proving that up to a subsequence,~$\{\mathcal{A}_{i}(x,\nabla v_{i})\}_{i\in\mathbb{N}}$ converges to~$\mathcal{A}(x,\nabla v)$ a.e. in~$\omega_{k}$. Take a function~$0\leq\varphi_{k}\in C^{\infty}_{c}(\omega_{k+1})\leq 1$ such that~$\varphi_{k}|_{\omega_{k}}\equiv 1$. For all sufficiently large~$i\in\mathbb{N}$, take~$\varphi_{k}(v_{i}-v)\in W^{1,p}_{c}(\omega_{k+1})$ as a test function and hence
\bals
&\int_{\omega_{k+1}}\varphi_{k}\mathcal{A}_{i}(x,\nabla v_{i})\cdot \nabla (v_{i}-v)\dx\\
&=-\int_{\omega_{k+1}}(v_{i}-v)\mathcal{A}_{i}(x,\nabla v_{i})\cdot \nabla\varphi_{k}\dx-\int_{\omega_{k+1}}V_{i}v_{i}^{p-1}\varphi_{k}(v_{i}-v)\dx.
\eals
For all sufficiently large~$i\in\mathbb{N}$, 
\bals
&\left|-\int_{\omega_{k+1}}(v_{i}-v)\mathcal{A}_{i}(x,\nabla v_{i})\cdot \nabla\varphi_{k}\dx\right|\leq\beta_{\omega_{k+1}}C(\varphi_{k})\int_{\omega_{k+1}}|\nabla v_{i}|^{p-1}|v_{i}-v|\dx\\
&\leq\beta_{\omega_{k+1}}C(\varphi_{k})\left(\int_{\omega_{k+1}}|v_{i}-v|^{p}\dx\right)^{1/p}\left(\int_{\omega_{k+1}}|\nabla v_{i}|^{p}\dx\right)^{1-1/p}\xrightarrow[i\longrightarrow\infty]{\mbox{up to a subsequence}}0,
\eals
where the convergence is because~$\{v_{i}\}_{i\in\mathbb{N}}$ converges uniformly in~$\omega_{k+1}$ to~$v$ up to a subsequence and~$v_{i}$ is bounded in~$W^{1,p}(\omega_{k+1})$ with respect to all sufficiently large positive integers~$i$. Because~$\Vert V_{i}\Vert_{M^{q}(p;\,\omega_{k+2})}$ is bounded for all sufficiently large~$i\in\mathbb{N}$,~$\int_{\omega_{k+1}}|V_{i}|\dx$ is bounded for all sufficiently large~$i\in\mathbb{N}$, which is concluded by a finite covering argument for~$p\leq n$ and by definition for~$p>n$. In~$\omega_{k+1}$, recall that~$\varphi_{k}$ is bounded, that~$\{v_{i}\}_{i\in\mathbb{N}}$ converges uniformly to~$v$ up to a subsequence, and that~$v_{i}$ is uniformly bounded for all sufficiently large~$i\in\mathbb{N}$. Then up to a subsequence,
\bals
\lim_{i\rightarrow\infty}\int_{\omega_{k+1}}V_{i}v_{i}^{p-1}\varphi_{k}(v_{i}-v)\dx=0.
\eals
Immediately, we conclude that up to a subsequence,
\bal\label{f1}
\lim_{i\rightarrow\infty}\int_{\omega_{k+1}}\varphi_{k}\mathcal{A}_{i}(x,\nabla v_{i})\cdot \nabla (v_{i}-v)\dx=0.
\eal
Moreover, up to a subsequence, since~$\{v_{i}\}_{i\in\mathbb{N}}$ converges weakly in~$W^{1,p}(\omega_{k+1})$ to~$v$, 
\bal\label{f2}
\lim_{i\rightarrow\infty}\int_{\omega_{k+1}}\varphi_{k}\mathcal{A}(x,\nabla v)\cdot(\nabla v_{i}-\nabla v)\dx=0.
\eal
For all sufficiently large~$i\in\mathbb{N}$, by H\"older's inequality and Minkowski's inequality,
\bal\label{f3p}
&\left|\int_{\omega_{k+1}}\varphi_{k}(x)(\mathcal{A}_{i}(x,\nabla v)-\mathcal{A}(x,\nabla v))\cdot(\nabla v_{i}-\nabla v)\dx\right|\notag\\
&\leq\int_{\omega_{k+1}}|\mathcal{A}_{i}(x,\nabla v)-\mathcal{A}(x,\nabla v)||\nabla v_{i}-\nabla v|\dx\notag\\
&\leq\left(\int_{\omega_{k+1}}|\mathcal{A}_{i}(x,\nabla v)-\mathcal{A}(x,\nabla v)|^{p'}\dx\right)^{1/p'}\left(\int_{\omega_{k+1}}|\nabla v_{i}-\nabla v|^{p}\dx\right)^{1/p}\notag\\
&\leq \left(\int_{\omega_{k+1}}|\mathcal{A}_{i}(x,\nabla v)-\mathcal{A}(x,\nabla v)|^{p'}\dx\right)^{1/p'}(\Vert\nabla v_{i}\Vert_{L^{p}(\omega_{k+1})}+\Vert\nabla v\Vert_{L^{p}(\omega_{k+1})}).
\eal
By assumption, for almost all~$x\in\omega_{k+1}$,~$\mathcal{A}_{i}(x,\nabla v)$ converges to~$\mathcal{A}(x,\nabla v)$ as~$i\rightarrow\infty$ and for all sufficiently large~$i\in\mathbb{N}$,
\bals
&|\mathcal{A}_{i}(x,\nabla v)-\mathcal{A}(x,\nabla v)|^{p'}\leq(|\mathcal{A}_{i}(x,\nabla v)|+|\mathcal{A}(x,\nabla v)|)^{p'}\\
&\leq 2^{p'-1}|\mathcal{A}_{i}(x,\nabla v)|^{p'}+2^{p'-1}|\mathcal{A}(x,\nabla v)|^{p'}\leq C|\nabla v|^{p}.
\eals
By the dominated convergence theorem,
\bals
&\lim_{i\rightarrow\infty}\left(\int_{\omega_{k+1}}|\mathcal{A}_{i}(x,\nabla v)-\mathcal{A}(x,\nabla v)|^{p'}\dx\right)^{1/p'}=0.
\eals
Recall that with respect to all sufficiently large positive integers~$i$,~$v_{i}$ is bounded in $W^{1,p}(\omega_{k+1})$.
Then \eqref{f3p} leads to
\bal\label{f3}
&\lim_{i\rightarrow\infty}\int_{\omega_{k+1}}\varphi_{k}(x)(\mathcal{A}_{i}(x,\nabla v)-\mathcal{A}(x,\nabla v))\cdot(\nabla v_{i}-\nabla v)\dx=0.
\eal
For almost all~$x\in\omega_{k+1}$ and all sufficiently large~$i\in\mathbb{N}$, let
\bals
&E_{i}(x)\triangleq (\mathcal{A}_{i}(x,\nabla v_{i})-\mathcal{A}_{i}(x,\nabla v))\cdot(\nabla v_{i}-\nabla v),\\
&G_{i}(x)\triangleq(\mathcal{A}_{i}(x,\nabla v_{i})-\mathcal{A}(x,\nabla v))\cdot(\nabla v_{i}-\nabla v),\\
\mbox{and}\quad&I_{i}(x)\triangleq (\mathcal{A}_{i}(x,\nabla v)-\mathcal{A}(x,\nabla v))\cdot(\nabla v_{i}-\nabla v).
\eals For almost all~$x\in\omega_{k+1}$ and all sufficiently large~$i\in\mathbb{N}$,
\bals
E_{i}(x)&=(\mathcal{A}_{i}(x,\nabla v_{i})-\mathcal{A}(x,\nabla v))\cdot(\nabla v_{i}-\nabla v)-(\mathcal{A}_{i}(x,\nabla v)-\mathcal{A}(x,\nabla v))\cdot(\nabla v_{i}-\nabla v)\\
&=G_{i}(x)-I_{i}(x),
\eals
and by the strict monotonicity of~$\mathcal{A}_{i}$,~$E_{i}(x)\geq 0.$
Consequently, for all sufficiently large~$i\in\mathbb{N}$,
\bals
&0\leq\int_{\omega_{k}}E_{i}(x)\dx=\int_{\omega_{k}}(G_{i}(x)-I_{i}(x))\dx\leq\int_{\omega_{k+1}}\varphi_{k}(x)(G_{i}(x)-I_{i}(x))\dx\\
&=\int_{\omega_{k+1}}\varphi_{k}(x)G_{i}(x)\dx-\int_{\omega_{k+1}}\varphi_{k}(x)I_{i}(x)\dx\\
&=\int_{\omega_{k+1}}\varphi_{k}\mathcal{A}_{i}(x,\nabla v_{i})\cdot(\nabla v_{i}-\nabla v)\dx-\int_{\omega_{k+1}}\varphi_{k}\mathcal{A}(x,\nabla v)\cdot(\nabla v_{i}-\nabla v)\dx\\
&\quad-\int_{\omega_{k+1}}\varphi_{k}(x)I_{i}(x)\dx\xrightarrow[i\longrightarrow\infty]{\mbox{up to a subsequence}}0,
\eals
where the last step is due to \eqref{f1}, \eqref{f2}, and \eqref{f3}. The remaining proof follows Maz'ya's argument (see \cite[Lemma~3.73]{HKM} and \cite[Lemma~1]{Maz'ya76}). Now it is easy to see that there exists a Lebesgue null set~$Z$ such that up to a subsequence,~$E_{i}(x)$ converges to~$0$ as~$i\rightarrow\infty$ for all~$x\in\omega_{k}\setminus Z$. Recall that up to a larger null set,~$\mathcal{A}_{i}(x,\cdot)$ is pointwise equicontinuous in~$\R^{n}$ for every~$x\in\omega_{k}\setminus Z$ and all sufficiently large~$i\in\mathbb{N}$. Since~$v\in W^{1,p}(\omega_{k})$, up to a larger null set, we may also assume that for all~$x\in\omega_{k}\setminus Z$,~$|\nabla v(x)|<\infty$. Next we fix an arbitrary~$x\in\omega_{k}\setminus Z$ and denote the two constants in the local uniform ellipticity and boundedness conditions of~$\mathcal{A}_{i}$ in~$\omega_{k}$ by~$\bar{\alpha}_{\omega_{k}}$ and~$\bar{\beta}_{\omega_{k}}$ for all sufficiently large~$i\in\mathbb{N}$, respectively. For all sufficiently large~$i\in\mathbb{N}$,
\bals
E_{i}(x)&=(\mathcal{A}_{i}(x,\nabla v_{i})-\mathcal{A}_{i}(x,\nabla v))\cdot(\nabla v_{i}-\nabla v)\\
&\geq \mathcal{A}_{i}(x,\nabla v_{i})\cdot\nabla v_{i}+\mathcal{A}_{i}(x,\nabla v)\cdot\nabla v-\bar{\beta}_{\omega_{k}}(|\nabla v_{i}|^{p-1}|\nabla v|+|\nabla v|^{p-1}|\nabla v_{i}|)\\
&\geq \bar{\alpha}_{\omega_{k}}|\nabla v_{i}|^{p}-\bar{\beta}_{\omega_{k}}(|\nabla v_{i}|^{p-1}|\nabla v|+|\nabla v|^{p-1}|\nabla v_{i}|)\\
&\geq \bar{\alpha}_{\omega_{k}}|\nabla v_{i}|^{p}-\bar{\beta}_{\omega_{k}}(|\nabla v|+|\nabla v|^{p-1}|)(|\nabla v_{i}|^{p-1}+|\nabla v_{i}|).
\eals
Therefore, up to a subsequence, which is independent of~$x\in\omega_{k}\setminus Z$,~$\{\nabla v_{i}(x)\}_{i\in\mathbb{N}}$ is bounded. Up to a subsequence, which is independent of~$x\in\omega_{k}\setminus Z$, let~$\eta\in\R^{n}$ be an arbitrary limit point of~$\{\nabla v_{i}(x)\}_{i\in\mathbb{N}}$. Then by the pointwise equicontinuity of~$\mathcal{A}_{i}(x,\cdot)$ for all sufficiently large~$i\in\mathbb{N}$, we obtain, up to a subsequence,\bal\label{aicon}
\lim_{i\rightarrow\infty}\mathcal{A}_{i}(x,\nabla v_{i})=\mathcal{A}(x,\eta).
\eal
Furthermore, up to a subsequence, which is independent of~$x\in\omega_{k}\setminus Z$, every limit point~$\eta\in\R^{n}$ of~$\{\nabla v_{i}(x)\}_{i\in\mathbb{N}}$ is~$\nabla v(x)$, which is because, up to a subsequence,
\bals
&(\mathcal{A}(x,\eta)-\mathcal{A}(x,\nabla v))\cdot(\eta-\nabla v)\\
&=\lim_{i\rightarrow\infty}(\mathcal{A}_{i}(x,\nabla v_{i})-\mathcal{A}_{i}(x,\nabla v))\cdot(\nabla v_{i}-\nabla v)=\lim_{i\rightarrow\infty}E_{i}(x)=0.
\eals
It follows that up to a subsequence,~$\lim_{i\rightarrow\infty}\nabla v_{i}=\nabla v$ a.e. in~$\omega_{k}$. Combining this and \eqref{aicon}, we conclude that up to a subsequence,~$\lim_{i\rightarrow\infty}\mathcal{A}_{i}(x,\nabla v_{i})=\mathcal{A}(x,\nabla v)$ a.e. in~$\omega_{k}$.

Since~$v$ is nonnegative in~$\Omega$ and~$v(x_{0})=1$, by the local Harnack inequality \cite[Theorem 3.1]{HPR},~$v$ is positive in~$\omega_{k}$ for all~$k\in\mathbb{N}$ and hence in~$\Omega$.
\eproof
\section{Fuchsian-type singularity}
In this section, we study a special class $\mathbf{M}_{\zeta,\hat{R}}$ of positive solutions of minimal growth near $S_{\hat{R}}(0)$. By virtue of the Fuchsian-type singularity with respect to a norm $N\in C^{1}(\R^{n}\setminus\{0\})$ on $\R^{n}$, we prove a uniform Harnack inequality and investigate asymptotic behaviors of quotients of positive solutions of $Q'[w]=0$ in~$\Gamma_{\zeta,\hat{R}}$.

Recall that $\widehat{\R^{n}}$ is the one-point compactification of~$\R^{n}$ with the ideal point $\infty$, that $\widehat{\partial}\Omega$ is the boundary of~$\Omega$ in $\widehat{\R^{n}}$, and that~$\zeta\in\{0,\infty\}$.

\textbf{Recall that~$\zeta$ is an isolated point of~$\widehat{\partial}\Omega$.} Recall that for all~$r,R>0$,~$\mathbb{A}_{r,R}=\{x\in\R^{n}~|~r<|x|<R\}$,$$I_{\zeta,R}=\begin{cases}
(0,R)& \mbox{if}~\zeta=0,\\
(R,\infty)& \mbox{if}~\zeta=\infty,	
\end{cases}
\quad\mbox{and}\quad \Gamma_{\zeta,R}=\begin{cases}
B_{R}(0)\setminus\{0\}& \mbox{if}~\zeta=0,\\
\R^{n}\setminus D_{R}(0)& \mbox{if}~\zeta=\infty.	
\end{cases}
$$
Recall that we fix~$\hat{R}>0$ such that~$\Gamma_{\zeta,R}\subseteq\Omega$ for all~$R\in I_{\zeta,\hat{R}}$ and that~$\Gamma_{\zeta,\hat{R}}\subseteq\Omega$.

Now we define positive solutions of minimal growth around $S_{\hat{R}}(0)$. 
  \bdefinition
  \begin{itemize}
 \item[\emph{(1)}]\emph{Let~$\zeta=0$. A positive solution~$u$ of~$Q'[w]=0$ in $B_{\hat{R}}(0)\sms\{0\}$ has \emph{minimal growth near $S_{\hat{R}}(0)$} if for every admissible compact subset $K$ of $B_{\hat{R}}(0)$ with $0\in\mathring{K}$ and for all positive solutions~$v\in C\left(B_{\hat{R}}(0)\setminus \mathring{K}\right)$ of~$Q'[w]=g$ in $B_{\hat{R}}(0)\setminus K$,$$u\leq v\quad\mbox{on}\quad\partial K\Longrightarrow u\leq v\quad\mbox{in}\quad B_{\hat{R}}(0)\setminus K,$$ where\begin{align*}
 g\in\begin{cases}
M^{q}_{\loc}(p;B_{\hat{R}}(0)\sms\{0\})&\mbox{if}~p\neq n,\\
L^{\bar{\rho}}_{\loc}(B_{\hat{R}}(0)\sms\{0\})~(\bar{\rho}>1)&\mbox{if}~p=n,
\end{cases}
\end{align*}
and $g$ is nonnegative a.e. in $B_{\hat{R}}(0)\setminus K$. We denote by $\mathbf{M}_{0,\hat{R}}$ the set of all such positive solutions.}
\item[\emph{(2)}]\emph{Let~$\zeta=\infty$. A positive solution~$u$ of~$Q'[w]=0$ in~$\R^{n}\setminus D_{\hat{R}}(0)$ has \emph{minimal growth near $S_{\hat{R}}(0)$} if for every bounded Lipschitz domain~$\omega\supseteq D_{\hat{R}}(0)$ and for all positive solutions~$v\in C\left(\overline{\omega}\setminus D_{\hat{R}}(0)\right)$ of~$Q'[w]=g$  in $\omega\setminus D_{\hat{R}}(0)$,$$u\leq v\quad\mbox{on}\quad\partial\omega\Longrightarrow u\leq v\quad\mbox{in}\quad \omega\setminus D_{\hat{R}}(0),$$ where\begin{align*}
 g\in\begin{cases}
M^{q}(p;\omega\setminus D_{\hat{R}}(0))&\mbox{if}~p\neq n,\\
L^{\bar{\rho}}\left(\omega\setminus D_{\hat{R}}(0)\right)~(\bar{\rho}>1)&\mbox{if}~p=n,
\end{cases}
\end{align*}
and $g$ is nonnegative a.e. in~$\omega\setminus D_{\hat{R}}(0)$. We denote by $\mathbf{M}_{\infty,\hat{R}}$ the set of all such positive solutions.}
\end{itemize}
  \edefinition
Under some further assumptions, we show that $\mathbf{M}_{\zeta,\hat{R}}\neq\emptyset$. 
  \blemma\label{bmgc}
 Suppose that \Cref{C1alpha} holds in $\Gamma_{\zeta,\hat{R}}$, that~$q>n$ for~$p<2$, that~$V\in\mathbb{M}^{q}_{\loc}(p;\Gamma_{\zeta,\hat{R}})$, and that $Q$ is nonnegative in $\Gamma_{\zeta,\hat{R}}$. Then~$\mathbf{M}_{\zeta,\hat{R}}\neq\emptyset$.
  \elemma
  \bproof
  When~$\zeta=0$, the proof is the same as that of \cite[Theorem 7.4]{HPR} and hence omitted. 

Let $\zeta=\infty$. The proof is similar to those of \cite[Proposition~3.9]{Giri1} and \cite[Theorem~9.7]{Hou}.

Consider two sequences $\{r_{k}\}_{k\in\N}$ decreasing to $\hat{R}$ and $\{R_{k}\}_{k\in\N}$ increasing to~$\infty$ such that $r_{k}<R_{k}$ for all $k\in\mathbb{N}$. Take a sequence $\{x_{k}\}_{k\in\mathbb{N}}\subseteq\R^{n}\setminus D_{\hat{R}}(0)$ converging to~$\infty$ with~$x_{1}\in\mathbb{A}_{r_{1},R_{1}}$ and~$x_{k}\in\mathbb{A}_{r_{k},R_{k}}\setminus\overline{\mathbb{A}_{r_{k-1},R_{k-1}}}$ for all~$k\in\N$ ($k\geq 2$). Note that for some sequence~$\{\delta_{k}\}_{k\in\mathbb{N}}\subseteq(0,1)$,~$D_{\delta_{k}}(x_{k})\subseteq\mathbb{A}_{r_{k},R_{k}}$ for every $k\in\N$ and~$\left\{\mathbb{A}_{r_{k},R_{k}}\setminus D_{\delta_{k}}(x_{k})\right\}_{k\in\mathbb{N}}$ is still a Lipschitz exhaustion of $\R^{n}\setminus D_{\hat{R}}(0)$. For every~$k\in\mathbb{N}$, take a nonnegative~$f_{k}\in C^{\infty}_{c}(B_{\delta_{k}}(x_{k}))\setminus\{0\}$. Since~$Q$ is nonnegative in $\R^{n}\setminus D_{\hat{R}}(0)$, by \cite[Corollary 6.15]{HPR}, $\lambda_{1}(Q;\mathbb{A}_{r_{k},R_{k}})>0$ for all~$k\in\mathbb{N}$. Furthermore, by \cite[Theorem 4.22]{HPR}, for every~$k\in\mathbb{N}$, there exists a postive solution~$u_{k}\in W^{1,p}_{0}(\mathbb{A}_{r_{k},R_{k}})$ of the equation $Q'[w]=f_{k}$ in~$\mathbb{A}_{r_{k},R_{k}}$. By \cite[Remark 5.5]{Hou3}, we get~$u_{k}\in C\left(\overline{\mathbb{A}_{r_{k},R_{k}}}\right)$ for all~$k\in\mathbb{N}$. Take~$x_{0}\in\mathbb{A}_{r_{1},R_{1}}\setminus D_{\delta_{1}}(x_{1})$ and let~$v_{k}\triangleq u_{k}/u_{k}(x_{0})$ for every~$k\in\mathbb{N}$. Since~$Q'[v_{k}]=0$ in~$\mathbb{A}_{r_{k},R_{k}}\setminus D_{\delta_{k}}(x_{k})$ for every~$k\in\mathbb{N}$, by the Harnack convergence principle \cite[Theorem 3.5]{HPR}, up to a subsequence,~$\{v_{k}\}_{k\in\mathbb{N}}$ converges locally uniformly in~$\R^{n}\setminus D_{\hat{R}}(0)$ to a positive solution~$v$ of the equation~$Q'[w]=0$.
 
Now we show that~$v\in\mathbf{M}_{\infty,\hat{R}}$. Fix an arbitrary bounded Lipschitz domain $\omega\supseteq D_{\hat{R}}(0)$ and an arbitrary positive solution~$\nu\in C\left(\overline{\omega}\setminus D_{\hat{R}}(0)\right)$ of~$Q'[w]=g$ in $\omega\setminus D_{\hat{R}}(0)$ such that $v\leq\nu$ on~$\partial\omega$, where\begin{align*}
 g\in\begin{cases}
M^{q}(p;\omega\setminus D_{\hat{R}}(0))&\mbox{if}~p\neq n,\\
L^{\bar{\rho}}\left(\omega\setminus D_{\hat{R}}(0)\right)~(\bar{\rho}>1)&\mbox{if}~p=n,
\end{cases}
\end{align*}
and $g$ is nonnegative a.e. in~$\omega\setminus D_{\hat{R}}(0)$. Fix an arbitrary~$\varepsilon>0$. Up to a subsequence, for all sufficiently large positive integers~$k$,~$v_{k}\leq (1+\varepsilon)\nu$ on~$\partial\omega\cup S_{r_{k}}(0)$. Note that~$\overline{\omega}\cap D_{\delta_{k}}(x_{k})=\emptyset$ for all sufficiently large~$k\in\mathbb{N}$. Up to a subsequence, for all sufficiently large positive integers~$k$,
\begin{align*}
\begin{cases}
		Q'[v_{k}]=0\leq (1+\varepsilon)^{p-1}g=Q'[(1+\varepsilon)\nu]&\mbox{in}~\omega\sms D_{r_{k}}(0),\\
			v_{k}\leq (1+\varepsilon)\nu&\mbox{on}~\partial\left(\omega\sms D_{r_{k}}(0)\right)=\partial\omega\cup S_{r_{k}}(0).
		\end{cases}
\end{align*}
By \cite[Corollary 6.15]{HPR},~$\lambda_{1}(Q;\omega\sms D_{r_{k}}(0))>0$ for all sufficiently large positive integers~$k$. By the weak comparison principle \cite[Theorem 4.25]{HPR} and \cite[Remark 5.5]{Hou3}, up to a subsequence, for all sufficiently large positive integers~$k$,~$v_{k}\leq (1+\varepsilon)\nu$ in~$\omega\sms D_{r_{k}}(0)$. Furthermore, $v\leq (1+\varepsilon)\nu$ in~$\omega\sms D_{\hat{R}}(0)$. Then~$v\leq\nu$ in~$\omega\sms D_{\hat{R}}(0)$. Therefore, $v\in\mathbf{M}_{\infty,\hat{R}}$.
  \eproof  
\bremark\label{vrrem}
\begin{itemize}
\item[\textup{(1)}]\emph{Let~$R>0$. Clearly,~$\Omega/R=\{x/R\in\R^{n}~|~x\in \Omega\}$ is still a domain. For almost all $x\in\Omega/R$, let $H^{(R)}(x,\cdot)\triangleq H(Rx,\cdot)$, let $\mathcal{A}^{(R)}(x,\cdot)\triangleq\mathcal{A}(Rx,\cdot)$, and let $V^{(R)}(x)\triangleq R^{p}V(Rx)$. Suppose that~$Q'$ has a Fuchsian-type singularity at~$\zeta$ with respect to a norm $N\in C^{1}(\R^{n}\setminus\{0\})$ on $\R^{n}$. Then for all~$R>0$, $\mathbb{A}_{a_{N}R,b_{N}R,N}/R=\mathbb{A}_{a_{N},b_{N},N}$ and for some~$R_{N}>0$ and all $R\in I_{\zeta,R_{N}}$, $\mathbb{A}_{a_{N}R,b_{N}R,N}\Subset\Gamma_{\zeta,\hat{R}}$, $\mathbb{A}_{a_{N},b_{N},N}\Subset\Gamma_{\zeta,\hat{R}}/R$, and the function,
\bals
\|V^{(R)}\|_{M^{q}(p;\mathbb{A}_{a_{N},b_{N},N})}=\begin{cases}
R^{p-n/q}\|V\|_{M^{q}(p;\mathbb{A}_{a_{N}R,b_{N}R,N})}&\mbox{if}~p\neq n,\\
\|V\|_{M^{q}(n;\mathbb{A}_{a_{N}R,b_{N}R,N})}&\mbox{if}~p=n,
\end{cases}
\eals
is bounded on~$I_{\zeta,R_{N}}$.}
\item[\textup{(2)}]\emph{Let~$\{R_{k}\}_{k\in\mathbb{N}}\subseteq (0,\infty)$ be a sequence converging to~$|\zeta|$ and suppose that~$0\notin\Omega$ if~$\zeta=\infty$. 
Then
\begin{align*}
\cup_{k\in\mathbb{N}}\Omega/R_{k}=\R^{n}\setminus\{0\}.
\end{align*} and $H^{(R_{k})}$ satisfies Assumptions~\ref{ass9} in~$\Omega/R_{k}$ for every~$k\in\mathbb{N}$.}
\end{itemize}
\eremark
\btheorem[{Uniform Harnack inequality}]\label{uhi}
Let~$u,v\in\mathcal{G}_{\zeta,\hat{R}}$. Suppose that~$Q'$ has a Fuchsian-type singularity at~$\zeta$ with respect to a norm $N\in C^{1}(\R^{n}\setminus\{0\})$ on~$\R^{n}$. Let $a_{N}<a_{N}'<1<b_{N}'<b_{N}$. Then there exist two positive constants~$R_{N}$ and~$C$ such that for all~$R\in I_{\zeta,R_{N}}$,
$$\sup_{\mathbb{A}_{a_{N}'R,b_{N}'R,N}}\frac{u}{v}\leq C\inf_{\mathbb{A}_{a_{N}'R,b_{N}'R,N}}\frac{u}{v},$$
where $C$ depends only on $n,p,q,\mathbb{A}_{a_{N}',b_{N}',N},\mathbb{A}_{a_{N},b_{N},N},\alpha_{\zeta,\hat{R}},\beta_{\zeta,\hat{R}}$, and the bound of\\$\|V^{(R)}\|_{M^{q}(p;\mathbb{A}_{a_{N},b_{N},N})}$ on $I_{\zeta,R_{N}}$.
\etheorem
\bproof
The proof is similar to that of \cite[Theorem 3.13]{Giri1}. For every $R>0$ and all $x\in\Gamma_{\zeta,\hat{R}}/R$, let~$u_{R}(x)=u(Rx)$ and let~$v_{R}(x)=v(Rx)$. Then $u_{R}$ and $v_{R}$ are positive solutions of the equation
\bals
-\dive\mathcal{A}^{(R)}(x,\nabla w)+V^{(R)}(x)w^{<p-1>}=0,
\eals
in $\Gamma_{\zeta,\hat{R}}/R$. By \Cref{vrrem} (1), there exists a positive constant~$R_{N}$ such that for all~$R\in I_{\zeta,R_{N}}$, $\mathbb{A}_{a_{N}R,b_{N}R,N}\Subset\Gamma_{\zeta,\hat{R}}$, $\mathbb{A}_{a_{N},b_{N},N}\Subset\Gamma_{\zeta,\hat{R}}/R$, and $\|V^{(R)}\|_{M^{q}(p;\mathbb{A}_{a_{N},b_{N},N})}$ is bounded on $I_{\zeta,R_{0}}$. By the local Harnack inequality \cite[Theorem 3.1]{HPR}, for all~$R\in I_{\zeta,R_{N}}$,
\bals\sup_{\mathbb{A}_{a_{N}'R,b_{N}'R,N}}\frac{u}{v}=\sup_{\mathbb{A}_{a_{N}',b_{N}',N}}\frac{u_{R}}{v_{R}}&\leq\frac{\sup_{\mathbb{A}_{a_{N}',b_{N}',N}}u_{R}}{\inf_{\mathbb{A}_{a_{N}',b_{N}',N}}v_{R}}\\
&\leq\frac{C\inf_{\mathbb{A}_{a_{N}',b_{N}',N}}u_{R}}{\sup_{\mathbb{A}_{a_{N}',b_{N}',N}}v_{R}} \leq C\inf_{\mathbb{A}_{a_{N}',b_{N}',N}}\frac{u_{R}}{v_{R}}= C\inf_{\mathbb{A}_{a_{N}'R,b_{N}'R,N}}\frac{u}{v},
\eals
where $C$ depends only on $n,p,q,\mathbb{A}_{a_{N}',b_{N}',N},\mathbb{A}_{a_{N},b_{N},N},\alpha_{\zeta,\hat{R}},\beta_{\zeta,\hat{R}}$, and the bound of\\$\|V^{(R)}\|_{M^{q}(p;\mathbb{A}_{a_{N},b_{N},N})}$ on $I_{\zeta,R_{N}}$.
\eproof
\bdefinition
\emph{Suppose that~$Q$ is nonegative in~$\Gamma_{\zeta,\hat{R}}$. Let~$N\in C^{1}(\R^{n}\setminus\{0\})$ be a norm on~$\R^{n}$ and let~$u,v\in\mathcal{G}_{\zeta,\hat{R}}$. For all~$R\in(0,\infty)$ such that~$S_{R,N}\subseteq\Gamma_{\zeta,\hat{R}}$, we define$$m_{R,N}\triangleq m_{R,N}(u,v)\triangleq\min_{S_{R,N}}\frac{u}{v}\quad\mbox{and}\quad M_{R,N}\triangleq M_{R,N}(u,v)\triangleq\max_{S_{R,N}}\frac{u}{v}.$$When~$N=|\cdot|$, for all~$R\in I_{\gz,\hat{R}}$, we omit the subscript `$|\cdot|$'.} 
\edefinition
Now we are in a position to consider asymptotic behaviors of quotients of positive solutions of $Q'[w]=0$ in~$\Gamma_{\zeta,\hat{R}}$.
\blemma\label{monolem}
Suppose that~$Q$ is nonegative in~$\Gamma_{\zeta,\hat{R}}$, that~$V\in\mathbb{M}^{q}_{\loc}(p;\Gamma_{\zeta,\hat{R}})$, that \Cref{C1alpha} holds in~$\Gamma_{\zeta,\hat{R}}$, that $q>n$ for~$p<2$, and that~$Q'$ has a Fuchsian-type singularity at~$\zeta$ with respect to a norm $N\in C^{1}(\R^{n}\setminus\{0\})$ on~$\R^{n}$. Let~$u,v\in\mathcal{G}_{\zeta,\hat{R}}$.
\begin{itemize}
\item [\emph{(i)}] The functions~$m_{R,N}$ and~$M_{R,N}$ of~$R$ are eventually monotone around~$|\zeta|$. In particular, both $\myd{\lim_{R\rightarrow|\zeta|}m_{R,N}}$ and~$\myd{\lim_{R\rightarrow|\zeta|}M_{R,N}}$ exist in the wide sense.

\item [\emph{(ii)}]Let~$m\triangleq\myd{\lim_{R\rightarrow|\zeta|}m_{R,N}}$, let~$M\triangleq\myd{\lim_{R\rightarrow|\zeta|}M_{R,N}}$, and let~$u,v\in\mathbf{M}_{\zeta,\hat{R}}$. As~$R$ goes to~$|\zeta|$,~$m_{R,N}$ decreases to~$m$ and~$M_{R,N}$ increases to~$M$ and~$0<m\leq M<\infty$. Moreover, for all~$x\in\Gamma_{\zeta,\hat{R}}$,$$mv(x)\leq u(x)\leq Mv(x).$$
\end{itemize}
\elemma
\bproof
The proof is similar to those of \cite[Lemma~4.2 and Corollary~4.1]{Fraas1}. See also \cite[Remark~3.15]{Giri1}.

(i): 
Clearly, for all~$R\in I_{\zeta,\hat{R}}$ such that~$S_{R,N}\subseteq\Gamma_{\zeta,\hat{R}}$, both $m_{R,N}$ and $M_{R,N}$ are positive.

Let~$\zeta=\infty$. Let~$\{R_{k}\}_{k\in\mathbb{N}}\subseteq\left(\hat{R},\infty\right)$ be any sequence increasing to~$\infty$ such that~$S_{R_{k},N}\subseteq\R^{n}\sms D_{\hat{R}}(0)$ for all~$k\in\mathbb{N}$. For every~$k,j\in\mathbb{N}$ and all~$x\in S_{R_{k},N}\cup S_{R_{k+j},N}$,
$$v(x)\min\{m_{R_{k},N},m_{R_{k+j},N}\}\leq u(x)\leq v(x)\max\{M_{R_{k},N},M_{R_{k+j},N}\}.$$
By \cite[Corollary 6.15]{HPR}, for every~$k,j\in\mathbb{N}$, the principal eigenvalue~$\lambda_{1}(Q;B_{R_{k+j},N}\setminus D_{R_{k},N})>0$.
By the weak comparison principle \cite[Theorem 4.25]{HPR} and \cite[Remark 5.5]{Hou3}, for every~$k,j\in\mathbb{N}$ and all~$x\in D_{R_{k+j},N}\setminus B_{R_{k},N}$,
$$v(x)\min\{m_{R_{k},N},m_{R_{k+j},N}\}\leq u(x)\leq v(x)\max\{M_{R_{k},N},M_{R_{k+j},N}\}.$$
Then for every~$k,i,j\in\mathbb{N}$ with~$i\leq j$ and all~$x\in S_{R_{k+i},N}$,
$$v(x)\min\{m_{R_{k},N},m_{R_{k+j},N}\}\leq u(x)\leq v(x)\max\{M_{R_{k},N},M_{R_{k+j},N}\},$$
and hence
$$m_{R_{k+i},N}\geq\min\{m_{R_{k},N},m_{R_{k+j},N}\}\quad\mbox{and}\quad M_{R_{k+i},N}\leq\max\{M_{R_{k},N},M_{R_{k+j},N}\}.$$
If~$\{m_{R_{k},N}\}_{k\in\mathbb{N}}$ is not increasing or~$\{M_{R_{k},N}\}_{k\in\mathbb{N}}$ is not decreasing, it can be easily checked that~$\{m_{R_{k},N}\}_{k\in\mathbb{N}}$ is eventually decreasing or~$\{M_{R_{k},N}\}_{k\in\mathbb{N}}$ is eventually increasing, respectively.
Therefore, both~$\{m_{R_{k},N}\}_{k\in\mathbb{N}}$ and~$\{M_{R_{k},N}\}_{k\in\mathbb{N}}$ are eventually monotone.  Then the functions~$m_{R,N}$ and~$M_{R,N}$ of~$R$ are eventually monotone around~$\infty$.

The case~$\zeta=0$ is similar and hence omitted.

(ii): Let~$\zeta=0$ and let~$u,v\in\mathbf{M}_{0,\hat{R}}$. Consider any positive numbers~$R<R'<\hat{R}$ such that $D_{R',N}\subseteq B_{\hat{R}}(0)$. Note that for all~$x\in S_{R,N}$,~$v(x)\leq\frac{u(x)}{m_{R,N}}$ and~$M_{R,N}v(x)\geq u(x)$. Since~$u,v\in\mathbf{M}_{0,\hat{R}}$, for all~$x\in S_{R',N}$,~$v(x)\leq\frac{u(x)}{m_{R,N}}$ and~$M_{R,N}v(x)\geq u(x)$. Therefore,~$m_{R,N}\leq m_{R',N}$ and~$M_{R,N}\geq M_{R',N}$. 

Let~$\zeta=\infty$ and let~$u,v\in\mathbf{M}_{\infty,\hat{R}}$. Consider any positive numbers~$\hat{R}<R<R'$ such that~$D_{\hat{R}}(0)\subseteq B_{R,N}$. Note that for all~$x\in S_{R',N}$,~$v(x)\leq\frac{u(x)}{m_{R',N}}$ and~$M_{R',N}v(x)\geq u(x)$. Since~$u,v\in\mathbf{M}_{\infty,\hat{R}}$, for all~$x\in S_{R,N}$,~$v(x)\leq\frac{u(x)}{m_{R',N}}$ and~$M_{R',N}v(x)\geq u(x)$. Therefore,~$m_{R',N}\leq m_{R,N}$ and~$M_{R',N}\geq M_{R,N}$. 

Let $a_{N}<a_{N}'<1<b_{N}'<b_{N}$. By \Cref{uhi}, there exist two positive constants~$R_{N}$ and~$C$ such that for all~$R\in I_{\zeta,R_{N}}$,
$$M_{R,N}\leq\sup_{\mathbb{A}_{a_{N}'R,b_{N}'R,N}}\frac{u}{v}\leq C\inf_{\mathbb{A}_{a_{N}'R,b_{N}'R,N}}\frac{u}{v}\leq Cm_{R,N}.$$ Then $M\leq Cm$. Clearly,~$m$ is finite and~$M$ is positive. Therefore,~$m$ is positive and~$M$ is finite.

Let~$N=|\cdot|$ be the Euclidean norm. The last result is immediate.
\eproof
\bdefinition
\emph{A subset~$E_{N}$ of $\Gamma_{\zeta,\hat{R}}$ is called \emph{essential} at~$\zeta$ with respect to a norm $N\in C^{1}(\R^{n}\setminus\{0\})$ on~$\R^{n}$ if there exist two real numbers~$\rho_{N}\in(0,1)$ and~$\sigma_{N}\in (1,\infty)$ and a sequence~$\{R_{k}\}_{k\in\mathbb{N}}\subseteq (0,\infty)$ converging to~$|\zeta|$ such that~$E_{N}=\cup_{k\in\mathbb{N}}\mathbb{A}_{\rho_{N} R_{k},\sigma_{N} R_{k},N}$, where~$\mathbb{A}_{\rho_{N} R_{k},\sigma_{N} R_{k},N}\Subset\Gamma_{\zeta,\hat{R}}$ for all~$k\in\N$.} 
\edefinition
\bdefinition
\emph{Let~$E_{N}$ be an essential subset of~$\Gamma_{\zeta,\hat{R}}$ at~$\zeta$ with respect to a norm $N\in C^{1}(\R^{n}\setminus\{0\})$ on~$\R^{n}$. The operator $Q'$ has a \emph{Fuchsian-type singularity} on $E_{N}$ if  \begin{itemize}
    \item  there exist two positive constants~$\alpha_{E_{N}}$ and~$\beta_{E_{N}}$ such that for almost all~$x\in E_{N}$ and all~$\xi\in\R^{n}$,
    $$\alpha_{E_{N}}|\xi|^{p}\leq\mathcal{A}(x,\xi)\cdot\xi\quad\mbox{and}\quad
						|\mathcal{A}(x,\xi)|\leq \beta_{E_{N}}\,|\xi|^{{p}-1};$$
     \item the following sequence $\{\Theta_{k,N}\}_{k\in\N}$ is bounded:\begin{align*}
\Theta_{k,N}\triangleq\begin{cases}
\|N^{p-n/q}V\|_{M^{q}(p;\mathbb{A}_{\rho_{N} R_{k},\sigma_{N} R_{k},N})}&\mbox{if}~p\neq n,\\
\|V\|_{M^{q}(n;\mathbb{A}_{\rho_{N} R_{k},\sigma_{N} R_{k},N})}&\mbox{if}~p=n.
\end{cases}
\end{align*}
\end{itemize}}
\edefinition
\bremark
\emph{Let~$E_{N}$ be an essential subset of~$\Gamma_{\zeta,\hat{R}}$ at~$\zeta$ with respect to a norm $N\in C^{1}(\R^{n}\setminus\{0\})$ on~$\R^{n}$. Suppose that the operator $Q'$ has a Fuchsian-type singularity on $E_{N}$. For all~$k\in\N$, $\mathbb{A}_{\rho_{N} R_{k},\sigma_{N} R_{k},N}/R_{k}=\mathbb{A}_{\rho_{N},\sigma_{N},N}$. Then $\mathbb{A}_{\rho_{N},\sigma_{N},N}\Subset\Gamma_{\zeta,\hat{R}}/R_{k}$ for every $k\in\N$ and the sequence
\bals
\|V^{(R_{k})}\|_{M^{q}(p;\mathbb{A}_{\rho_{N},\sigma_{N},N})}=\begin{cases}
R_{k}^{p-n/q}\|V\|_{M^{q}(p;\mathbb{A}_{\rho_{N}R_{k},\sigma_{N}R_{k},N})}&\mbox{if}~p\neq n,\\
\|V\|_{M^{q}(n;\mathbb{A}_{\rho_{N}R_{k},\sigma_{N}R_{k},N})}&\mbox{if}~p=n,
\end{cases}
\eals is bounded. It can be easily checked that \Cref{monolem}, \Cref{uniquemini}, and \Cref{ldregular} still hold if we assume a Fuchsian-type singularity only on $E_{N}$.}
\eremark
The following corollary is a Picard-type principle for $\zeta=0$ and a positive Liouville-type theorem for $\zeta=\infty$, respectively.
\bcorollary\label{uniquemini}Suppose that~$Q$ is nonegative in~$\Gamma_{\zeta,\hat{R}}$, that~$V\in\mathbb{M}^{q}_{\loc}(p;\Gamma_{\zeta,\hat{R}})$, that \Cref{C1alpha} holds in~$\Gamma_{\zeta,\hat{R}}$, that~$q>n$ for~$p<2$, and that~$Q'$ has a Fuchsian-type singularity at~$\zeta$ with respect to a norm $N\in C^{1}(\R^{n}\sms\{0\})$ on $\R^{n}$. If the point~$\zeta$ is regular with respect to~$Q'$, then~$\mathbf{M}_{\zeta,\hat{R}}$ has a unique function up to a positive multiplicative constant. 
\ecorollary
\bproof
The proof is similar to that of \cite[Proposition 3.17]{Giri1}. The existence is guaranteed by \Cref{bmgc}. We only need to show the uniqueness. Let~$u,v\in\mathbf{M}_{\zeta,\hat{R}}$. By \Cref{monolem}~(ii), for all~$x\in\Gamma_{\zeta,\hat{R}}$,$$0<m\leq\frac{u(x)}{v(x)}\leq M<\infty.$$ Since the point~$\zeta$ is regular with respect to~$Q'$,~$\myd{\lim_{x\rightarrow\zeta}u(x)/v(x)}$ exists. Then~$m=M$ and~$u=mv$ in~$\Gamma_{\zeta,\hat{R}}$.
\eproof
\section{Asymptotic behaviors of positive Finsler $p$-harmo-\\nic functions and their quotients}
In this section, for a fixed norm $H=H(\xi)$ ($\xi\in\R^{n}$), we consider the corresponding Finsler $p$-Laplace equation $-\dive\mathcal{A}(\nabla u)=0$, and study asymptotic behaviors of quotients of positive Finsler $p$-harmonic functions around $0$ for $p\leq n$ and asymptotic behaviors of positive Finsler $p$-harmonic functions around $\infty$.
\bdefinition
\emph{Suppose that $H$ is a fixed norm. For all~$\eta\in\R^{n}$, let
$$H_{0}(\eta)\triangleq\sup_{\xi\in\R^{n}\sms\{0\}}\frac{\eta\cdot\xi}{H(\xi)},$$ which is called the \emph{dual norm} of~$H$.}
\edefinition
We recall special Finsler $p$-harmonic functions in the punctured space for $p\neq n$. See also \cite[Lemma 4.2]{Giri1}.
\blemma[{\cite[Example~4.8~(iii)]{Hou3}}]
Suppose that $H$ is a fixed norm. Then $H_{0}^{(p-n)/(p-1)}$ is Finsler~$p$-harmonic in~$\R^{n}\sms\{0\}$ for~$p\neq n$.
\elemma
\blemma\label{ndim}
Suppose that $H$ is a fixed norm. Then $\log H_{0}$ is Finsler~$n$-harmonic in~$\R^{n}\sms\{0\}$.
\elemma
\bproof
The proof is similar to that of \cite[Lemma~4.2]{Giri1}. Recall that~$H_{0}\in C^{1}(\R^{n}\sms\{0\})$ (see \cite[Example~4.8~(iii)]{Hou3}). By \cite[p.~1139]{Jaros}, for all~$\eta\in\R^{n}\setminus\{0\}$ and all~$t\in\R\sms\{0\}$,
\begin{align*}
\nabla H(t\eta)=\sgn(t)\nabla H(\eta)\quad\mbox{and}\quad H(\nabla H_{0}(\eta))=1,
\end{align*}
and for all $\eta\in\R^{n}$,
\begin{align*}
H\left(H_{0}(\eta)\nabla H_{0}(\eta)\right)\nabla H\left(H_{0}(\eta)\nabla H_{0}(\eta)\right)=\eta\quad\mbox{and}\quad\eta\cdot\nabla H_{0}(\eta)=H_{0}(\eta),
\end{align*} where $H(0)\nabla H(0)$, $H_{0}(0)\nabla H_{0}(0)$, and~$0\cdot\nabla H_{0}(0)$ are all defined as $0$. Then for every~$\eta\in\R^{n}\sms\{0\}$,
$$\nabla\log H_{0}(\eta)=\frac{\nabla H_{0}(\eta)}{H_{0}(\eta)},$$
\bals
\mathcal{A}(\nabla\log H_{0}(\eta))&=H(\nabla\log H_{0}(\eta))^{n-1}\nabla H(\nabla\log H_{0}(\eta))\\
&=H\left(\frac{\nabla H_{0}(\eta)}{H_{0}(\eta)}\right)^{n-1}\nabla H\left(\frac{\nabla H_{0}(\eta)}{H_{0}(\eta)}\right)\\
&=\frac{1}{H_{0}(\eta)^{n-1}}H\left(\nabla H_{0}(\eta)\right)^{n-1}\nabla H\left(\nabla H_{0}(\eta)\right)\\
&=\frac{1}{H_{0}(\eta)^{n}}H\left(\nabla H_{0}(\eta)\right)^{n-2}H\left(H_{0}(\eta)\nabla H_{0}(\eta)\right)\nabla H\left(H_{0}(\eta)\nabla H_{0}(\eta)\right)\\
&=\frac{1}{H_{0}(\eta)^{n}}H\left(\nabla H_{0}(\eta)\right)^{n-2}\eta={H_{0}(\eta)^{-n}}\eta,
\eals
and
\bals
\dive\mathcal{A}(\nabla\log H_{0}(\eta))&=\sum_{i=1}^{n}\frac{\partial(H_{0}(\eta)^{-n}\eta_{i})}{\partial\eta_{i}}\\
&=\sum_{i=1}^{n}\left(H_{0}(\eta)^{-n}-nH_{0}(\eta)^{-n-1}\eta_{i}\frac{\partial H_{0}(\eta)}{\partial\eta_{i}}\right)\\
&=nH_{0}(\eta)^{-n}-nH_{0}(\eta)^{-n-1}\eta\cdot\nabla H_{0}(\eta)\\
&=nH_{0}(\eta)^{-n}-nH_{0}(\eta)^{-n-1}H_{0}(\eta)=0,
\eals
which implies that $\log H_{0}$ is Finsler~$n$-harmonic in~$\R^{n}\sms\{0\}$.
\eproof
\bdefinition
\emph{Suppose that $H$ is a fixed norm. We define
\bals
\mu_{p,H}\triangleq\begin{cases}
H_{0}^{(p-n)/(p-1)}&\mbox{if}~p\neq n,\\
-\log H_{0}&\mbox{if}~p=n.
\end{cases}
\eals}
\edefinition
The following theorem and \Cref{pncoro} are devoted to asymptotic behaviors of quotients of positive Finsler $p$-harmonic functions around $0$ for $p\leq n$.
\btheorem\label{asympn}
Let~$1<p\leq n$. Suppose that $H$ is a fixed norm, that \Cref{C1alphan} holds, and that~$u$ is a positive Finsler $p$-harmonic function in $B_{\hat{R}}(0)\sms\{0\}$ with a nonremovable singularity at~$0$. Then~$u(x)/\mu_{p,H}(x)$ converges to a positive constant as~$x\rightarrow 0$. 
\etheorem
\bproof
The proof is similar to that of \cite[Theorem 4.4]{Giri1}. By Serrin's result \cite[Theorem~1]{Serrin65},~$u\asymp\mu_{p,H}$ in~$B_{r_{0}}(0)\sms\{0\}$ for some positive number $r_{0}<\hat{R}$. Recall that for all~$r\in(0,\infty)$ such that~$S_{r,H_{0}}\subseteq B_{\hat{R}}(0)\sms\{0\}$,
$$M_{r,H_{0}}(u,\mu_{p,H})=\max_{S_{r,H_{0}}}\frac{u}{\mu_{p,H}}\quad\mbox{and}\quad m_{r,H_{0}}(u,\mu_{p,H})=\min_{S_{r,H_{0}}}\frac{u}{\mu_{p,H}}.$$
For all~$r\in(0,\infty)$ such that~$S_{r,H_{0}}\subseteq B_{r_{0}}(0)\sms\{0\}$, $M_{r,H_{0}}(u,\mu_{p,H})\asymp 1$ and~$m_{r,H_{0}}(u,\mu_{p,H})\asymp 1$.
By \Cref{monolem} (i), both $M=\lim_{r\rightarrow 0}M_{r,H_{0}}(u,\mu_{p,H})$ and $m=\lim_{r\rightarrow 0}m_{r,H_{0}}(u,\mu_{p,H})$ exist in~$(0,\infty)$. Take a sequence~$\{r_{k}\}_{k\in\N}\subseteq (0,r_{0})$ such that~$\{r_{k}\}_{k\in\N}$ strictly decreases to~$0$ as~$k\rightarrow\infty$ and that~$S_{r_{k},H_{0}}\subseteq B_{r_{0}}(0)\sms\{0\}$ for all~$k\in\N$. Take an arbitrary point~$x_{0}\in B_{1}(0)\sms\{0\}$. Then~$r_{k}x_{0}\in B_{r_{0}}(0)\sms\{0\}$ for all~$k\in\N$. For every~$k\in\N$ and all~$x\in (B_{r_{0}}(0)\sms\{0\})/r_{k}$, let
\bals
u_{k}(x)\triangleq\frac{u(r_{k}x)}{\mu_{p,H}(r_{k}x_{0})}\quad\mbox{and}\quad \mu_{p,H,k}(x)\triangleq\frac{\mu_{p,H}(r_{k}x)}{\mu_{p,H}(r_{k}x_{0})}=\begin{cases}
\frac{H_{0}(x)^{(p-n)/(p-1)}}{H_{0}(x_{0})^{(p-n)/(p-1)}}&\mbox{if}~p<n,\vspace{1mm}\\
\frac{-\log r_{k}-\log H_{0}(x)}{-\log r_{k}-\log H_{0}(x_{0})}&\mbox{if}~p=n.
\end{cases}
\eals
For every~$k\in\N$, both~$u_{k}$ and~$\mu_{p,H,k}$ are Finsler~$p$-harmonic in~$(B_{r_{0}}(0)\sms\{0\})/r_{k}$. For all~$k\in\N$,~$u_{k}(x_{0})\asymp 1$ and~$\mu_{p,H,k}(x_{0})=1$. By the Harnack convergence principle \cite[Theorem 3.5]{HPR}, there exists a strictly increasing sequence~$\{k_{j}\}_{j\in\N}$ of positive integers such that~$\{u_{k_{j}}\}_{j\in\N}$ and~$\{\mu_{p,H,k_{j}}\}_{j\in\N}$ converges locally uniformly in~$\R^{n}\sms\{0\}$ to two positive Finsler $p$-harmonic functions~$u_{\infty}$ and~$\mu_{p,H,\infty}$, respectively. It is a simple matter to verify that~$\{u_{k_{j}}/\mu_{p,H,k_{j}}\}_{j\in\N}$ converges locally uniformly in~$\R^{n}\sms\{0\}$ to~$u_{\infty}/\mu_{p,H,\infty}$. For every~$r\in(0,\infty)$, we may calculate:
\bals
&\max_{S_{r,H_{0}}}\frac{u_{\infty}}{\mu_{p,H,\infty}}=\lim_{j\rightarrow\infty}\max_{S_{r,H_{0}}}\frac{u_{k_{j}}}{\mu_{p,H,k_{j}}}=\lim_{j\rightarrow\infty}\max_{x\in S_{r,H_{0}}}\frac{u(r_{k_{j}}x)}{\mu_{p,H}(r_{k_{j}}x)}\\
&=\lim_{j\rightarrow\infty}\max_{x\in S_{rr_{k_{j}},H_{0}}}\frac{u(x)}{\mu_{p,H}(x)}=\lim_{j\rightarrow\infty}M_{rr_{k_{j}},H_{0}}(u,\mu_{p,H})=M,
\eals
and
\bals
&\min_{S_{r,H_{0}}}\frac{u_{\infty}}{\mu_{p,H,\infty}}=\lim_{j\rightarrow\infty}\min_{S_{r,H_{0}}}\frac{u_{k_{j}}}{\mu_{p,H,k_{j}}}=\lim_{j\rightarrow\infty}\min_{x\in S_{r,H_{0}}}\frac{u(r_{k_{j}}x)}{\mu_{p,H}(r_{k_{j}}x)}\\
&=\lim_{j\rightarrow\infty}\min_{x\in S_{rr_{k_{j}},H_{0}}}\frac{u(x)}{\mu_{p,H}(x)}=\lim_{j\rightarrow\infty}m_{rr_{k_{j}},H_{0}}(u,\mu_{p,H})=m.
\eals
Note that
\bals
\mu_{p,H,\infty}(x)=\begin{cases}
\frac{H_{0}(x)^{(p-n)/(p-1)}}{H_{0}(x_{0})^{(p-n)/(p-1)}}&\mbox{if}~p<n,\\
1&\mbox{if}~p=n.
\end{cases}
\eals
For all~$r\in(0,\infty)$, on~$S_{r,H_{0}}$,
\bals
m\mu_{p,H,\infty}\leq u_{\infty}\leq M\mu_{p,H,\infty}.
\eals
Therefore, for all~$x\in\R^{n}\sms\{0\}$,
\bals
m\mu_{p,H,\infty}(x)\leq u_{\infty}(x)\leq M\mu_{p,H,\infty}(x).
\eals
Recall that for all~$x\in\R^{n}\setminus\{0\}$,~$H(\nabla H_{0}(x))=1$ (see the proof of \Cref{ndim}) and hence~$\nabla H_{0}(x)\neq 0$. 
\begin{itemize}
\item Let~$p<n$. According to \cite[Remark 2.13]{Hou3}, both $u_{\infty}$ and~$\mu_{p,H,\infty}$ are continuously differentiable in~$\R^{n}\sms\{0\}$. For all~$x\in\R^{n}\sms\{0\}$,
\bals
\nabla\mu_{p,H,\infty}&=\frac{(p-n)H_{0}(x)^{(p-n)/(p-1)-1}}{(p-1)H_{0}(x_{0})^{(p-n)/(p-1)}}\nabla H_{0}(x)\\
&=\frac{(p-n)\nabla H_{0}(x)}{(p-1)H_{0}(x)^{(n-1)/(p-1)}H_{0}(x_{0})^{(p-n)/(p-1)}}\neq 0.
\eals
By \Cref{C1alphan}, both $D\mathcal{A}(\nabla\mu_{p,H,\infty})$ and $D\mathcal{A}(-\nabla\mu_{p,H,\infty})$ are positive definite in~$\R^{n}\sms\{0\}$.
By the tangency principle \cite[Theorem 2.5.2]{Serrin},
\bals
m\mu_{p,H,\infty}=u_{\infty}= M\mu_{p,H,\infty}.
\eals
As a result,~$m=M$.
\item Let~$p=n$. For all~$r\in(0,\infty)$, $M=\max_{S_{r,H_{0}}}u_{\infty}$ and $m=\min_{S_{r,H_{0}}}u_{\infty}$. In $\R^{n}\sms\{0\}$, $m\leq u_{\infty}\leq M$. By the strong maximum principle \cite[6.5]{HKM},~$m=u_{\infty}=M$.
\end{itemize}
Then$$\lim_{r\rightarrow0}\max_{S_{r,H_{0}}}\frac{u}{\mu_{p,H}}=M=m=\lim_{r\rightarrow0}\min_{S_{r,H_{0}}}\frac{u}{\mu_{p,H}}.$$
It follows that
\bals\lim_{x\rightarrow0}\frac{u(x)}{\mu_{p,H}(x)}=M=m.\qquad\qedhere\eals
\eproof
Since any two norms on~$\R^{n}$ are equivalent, the following lemma can be directly derived from \cite[Theorem 4.6]{Giri1} and \cite[Theorem 5.3]{HPR}. 
\blemma\label{cri_dim}
Suppose that $H$ is a fixed norm. If~$p\geq n$, then~$Q'_{0}$ is critical in~$\R^{n}$.
\elemma
Now we study asymptotic behaviors of positive Finsler $p$-harmonic functions around $\infty$.
\btheorem\label{theorem59}
Suppose that $H$ is a fixed norm, that \Cref{C1alphan} holds, and that~$u$ is a positive Finsler $p$-harmonic function in $\R^{n}\sms D_{\hat{R}}(0)$. Then~$\myd{\lim_{x\rightarrow\infty}u(x)}$ exists (possibly $\infty$). Furthermore,
\begin{itemize}
\item If~$p\geq n$, then $\myd{\lim_{x\rightarrow\infty}u(x)\neq 0}$.
\item If~$p<n$, then $\myd{\lim_{x\rightarrow\infty}u(x)\neq \infty}$.
\end{itemize}
\etheorem
\bproof
The proof is similar to those of \cite[Theorem 4.6 and Corollary 4.7]{Giri1}. By \Cref{monolem}, the functions~$m_{R}(u,1)=\min_{S_{R}(0)}u$ and~$M_{R}(u,1)=\max_{S_{R}(0)}u$ of~$R$ are eventually monotone near~$\infty$, and $m=\lim_{R\rightarrow\infty}m_{R}(u,1)$ and $M=\lim_{R\rightarrow\infty}M_{R}(u,1)$ exist. 

If~$m=\infty$, then $\lim_{x\rightarrow\infty}u(x)=\infty$. 

Suppose that~$m<\infty$. Fix an arbitrary~$\varepsilon>0$. There exists~$\bar{R}>\hat{R}$ such that~$u-m+\varepsilon$ is positive and clearly Finsler~$p$-harmonic in~$\R^{n}\sms D_{\bar{R}}(0)$. By \Cref{uhi}, there exists a positive constant~$C$ such that for all sufficiently large~$R$,
\bal\label{MRmR}
M_{R}(u,1)-m+\varepsilon&=M_{R}(u-m+\varepsilon,1)\notag\\
&\leq Cm_{R}(u-m+\varepsilon,1)=C(m_{R}(u,1)-m+\varepsilon).
\eal
Letting~$R\rightarrow\infty$ in \eqref{MRmR}, we see that~$0\leq M-m\leq (C-1)\varepsilon$. Since~$\varepsilon$ is an arbitrary positive constant, we get~$M=m<\infty$. Then~$\lim_{x\rightarrow\infty}u(x)=M=m<\infty$.

Let~$p\geq n$. By \Cref{cri_dim}, $Q'_{0}$ is critical in~$\R^{n}$. By \Cref{uniques},~$1$ is a unique (up to a positive multiplicative constant) positive supersolution of $Q'_{0}[u]=0$ in $\R^{n}$. By \cite[Theorems 6.9, 6.12, and 7.7]{HPR},~$1$ is a global minimal positive solution of~$Q'_{0}[v]=0$ in~$\R^{n}$. Then~$u\geq \min_{S_{\hat{R}+1}(0)}u>0$ in~$\R^{n}\sms D_{\hat{R}+1}(0)$. It follows that~$\lim_{x\rightarrow\infty}u(x)\neq 0$.

Next let~$p<n$ and suppose, contrary to our claim, that $\lim_{x\rightarrow\infty}u(x)=\infty$. 

For every~$k\in\N$, applying \cite[Theorem 4.24]{HPR} in~$B_{\hat{R}+k+1}(0)\sms D_{\hat{R}+1}(0)$ with any function $0\leq f=f_{k}\in C^{\infty}_{c}(\R^{n})\leq 1$ such that $f_{k}|_{S_{\hat{R}+1}(0)}=1$ and $f_{k}|_{S_{\hat{R}+k+1}(0)}=0$,~$\psi=0$, and~$\varphi=1$, we may find a nonnegative solution $v_{k}\in W^{1,p}(B_{\hat{R}+k+1}(0)\sms D_{\hat{R}+1}(0))\cap C(D_{\hat{R}+k+1}(0)\sms B_{\hat{R}+1}(0))$ to the Dirichlet problem
\begin{align*}
\begin{cases}		Q'_{0}[w]=0&~\mbox{in}~B_{\hat{R}+k+1}(0)\sms D_{\hat{R}+1}(0),\\		w=1&~\mbox{on}~S_{\hat{R}+1}(0),\\
w=0&~\mbox{on}~S_{\hat{R}+k+1}(0),
		\end{cases}
\end{align*}
such that $0\leq v_{k}\leq 1$ in $B_{\hat{R}+k+1}(0)\sms D_{\hat{R}+1}(0)$. By \cite[Theorem 3.1]{HPR}, for every~$k\in\N$,~$v_{k}$ is positive in $B_{\hat{R}+k+1}(0)\sms D_{\hat{R}+1}(0)$. By \cite[Corollary 6.15]{HPR},~$\lambda_{1}(Q_{0};B_{\hat{R}+k+1}(0)\sms D_{\hat{R}+1}(0))>0$ for every~$k\in\N$. By the weak comparison principle \cite[Theorem 4.25]{HPR} and \cite[Remark 5.5]{Hou3},~$\{v_{k}\}_{k\in\N}$ is increasing. By the Harnack convergence principle \cite[Theorem 3.5]{HPR}, in~$\R^{n}\sms D_{\hat{R}+1}(0)$, $\{v_{k}\}_{k\in\N}$ converges locally uniformly to a positive Finsler $p$-harmonic function~$v$. 

We claim that~$v\in\mathcal{M}_{\R^{n};D_{\hat{R}+1}(0)}$. Consider an arbitrary admissible compact subset $K$ of $\R^{n}$ with~$D_{\hat{R}+1}(0)\subseteq \mathring{K}$ and an arbitrary positive solution~$\tilde{v}\in C\left(\R^{n}\setminus \mathring{K}\right)$ of $Q'_{0}[w]=g$ in~$\R^{n}\setminus K$ such that $v\leq\tilde{v}$ on~$\partial K$, where~$g\in M^{q}_{\loc}(p;\R^{n})$ and $g$ is nonnegative a.e. in $\R^{n}\setminus K$. For all sufficiently large $k\in\mathbb{N}$,
\begin{align*}
\begin{cases}
		Q'_{0}[v_{k}]=0\leq g=Q'_{0}[\tilde{v}]&\mbox{in}~B_{\hat{R}+k+1}(0)\setminus K,\\		v_{k}\leq\tilde{v}&\mbox{on}~\partial(B_{\hat{R}+k+1}(0)\setminus K).
		\end{cases}
\end{align*}
By \cite[Corollary 6.15]{HPR}, for all sufficiently large~$k\in\mathbb{N}$, $\lambda_{1}(Q_{0};B_{\hat{R}+k+1}(0)\setminus K)>0$. By the weak comparison principle \cite[Theorem 4.25]{HPR} and \cite[Remark 5.5]{Hou3}, for all sufficiently large $k\in\mathbb{N}$, $v_{k}\leq \tilde{v}$ in $B_{\hat{R}+k+1}(0)\setminus K$. Then $v\leq\tilde{v}$ in $\R^{n}\setminus K$. 

Since $0\leq v_{k}\leq 1$ in $B_{\hat{R}+k+1}(0)\sms D_{\hat{R}+1}(0)$ for all~$k\in\N$ and $\{v_{k}\}_{k\in\N}$ increases to~$v$, we conclude that $v\leq 1$ in $\R^{n}\setminus D_{\hat{R}+1}(0)$. Fix any~$\varepsilon>0$. Because~$\lim_{x\rightarrow\infty}u(x)=\infty$, there exists~$k_{0}\in\N$ such that for all positive integers~$k\geq k_{0}$,~$1-\varepsilon u\leq v_{k}+\varepsilon$ on~$S_{\hat{R}+1}(0)\cup S_{\hat{R}+k+1}(0)$. By the weak comparison principle \cite[Theorem 4.25]{HPR} and \cite[Remark 5.5]{Hou3}, for all positive integers~$k\geq k_{0}$,~$1-\varepsilon u\leq v_{k}+\varepsilon$ in~$B_{\hat{R}+k+1}(0)\sms D_{\hat{R}+1}(0)$. It follows that~$1-\varepsilon u\leq v+\varepsilon$ in~$\R^{n}\sms D_{\hat{R}+1}(0)$. Letting~$\varepsilon\rightarrow 0$, we see that~$v\geq 1$ in~$\R^{n}\sms D_{\hat{R}+1}(0)$. Therefore,~$v=1$ in~$\R^{n}\sms D_{\hat{R}+1}(0)$. Furthermore,~$1\in\mathcal{M}_{\R^{n};D_{\hat{R}+1}(0)}$. Then~$\mu_{p,H}\geq\min_{S_{\hat{R}+2}}\mu_{p,H}>0$ in~$\R^{n}\sms D_{\hat{R}+2}(0)$, which is a contradiction to~$\lim_{x\rightarrow\infty}\mu_{p,H}=0$. In conclusion, when~$p<n$, $\lim_{x\rightarrow\infty}u(x)\neq\infty$. 
\eproof
\section{Limiting dilated equation and positive Liouville-type theorem}
In this section, we study a limiting dilated equation (operator) and for $p\leq n$, we prove a positive Liouville-type theorem. Under some further restrictions, we show that if~$0$ or~$\infty$ is regular with respect to the limiting dilated operator, then~$\zeta$ is regular with respect to~$Q'$. The positive Liouville-type theorem is based on a weak Fuchsian-type singularity.

We first define the limiting dilated equation (operator). For every~$R>0$, recall that for almost all~$x\in\Omega/R$, $\mathcal{A}^{(R)}(x,\cdot)\triangleq\mathcal{A}(Rx,\cdot)$, and~$V^{(R)}(x)=R^{p}V(Rx)$.
\bdefinition\label{limitingd}
\emph{Let~$\{R_{k}\}_{k\in\mathbb{N}}\subseteq (0,\infty)$ be a sequence converging to~$|\zeta|$ and suppose that~$0\notin\Omega$ if~$\zeta=\infty$. Suppose that for every domain compactly contained in $\R^{n}\setminus\{0\}$ and all sufficiently large~$k\in\mathbb{N}$, the two constants in the local uniform ellipticity and boundedness conditions of $\mathcal{A}^{(R_{k})}$ are independent of~$k\in\mathbb{N}$, that~$\mathcal{A}^{(R_{k})}(x,\cdot)$ is pointwise equicontinuous in~$\R^{n}$ for almost every~$x\in\R^{n}\setminus\{0\}$ and all sufficiently large~$k\in\mathbb{N}$, and that~$\lim_{k\rightarrow\infty}\mathcal{A}^{(R_{k})}(x,\xi)=\mathbf{A}(x,\xi)=\mathbf{A}_{\{R_{k}\}_{k\in\mathbb{N}}}(x,\xi)$ for almost all~$x\in\R^{n}\setminus\{0\}$ and all~$\xi\in\R^{n}$, where $\mathbf{A}$ satisfies all the conditions in \Cref{thm_1}. Suppose that $V^{(R_{k})}\in M^{q}_{\loc}(p;\Omega/R_{k})$ converges to~$\mathbf{V}=\mathbf{V}_{\{R_{k}\}_{k\in\mathbb{N}}}\in M^{q}_{\loc}(p;\R^{n}\setminus\{0\})$ weakly in~$M^{q}_{\loc}(p;\R^{n}\setminus\{0\})$ as~$k\rightarrow\infty$.}

\emph{We define the \emph{limiting dilated equation (operator)} with respect to~$Q'$ and~$\{R_{k}\}_{k\in\mathbb{N}}$ as
\begin{align*}
Q'_{\{R_{k}\}_{k\in\mathbb{N}}}[u]\triangleq-\dive\mathbf{A}(x,\nabla u)+\mathbf{V}u^{<p-1>}=0\quad\mbox{in}\quad\R^{n}\setminus\{0\}.
\end{align*}
}
\edefinition
\bexample
\emph{\begin{itemize}
    \item Let~$H(x,\xi)=\sqrt[p]{|\xi|_{s}^{p}+(|x|+1)^{p}|\xi|^{p}}$ ($2\leq s<\infty$) for all~$x\in\Omega$ and~$\xi\in\R^{n}$, where $|\xi|_{s}\triangleq\sqrt[s]{\sum_{i=1}^{n}|\xi_{i}|^{s}}$. Then $H$ satisfies \Cref{ass9}. Let $\zeta=0$ and let~$\{R_{k}\}_{k\in\mathbb{N}}\subseteq (0,\infty)$ be a sequence converging to~$0$. Then $\{\mathcal{A}^{(R_{k})}\}_{k\in\N}$ satisfies the conditions in \Cref{limitingd} with the limit operator $$\mathbf{A}(\xi)=|\xi|^{p-s}_{s}(|\xi_{1}|^{s-2}\xi_{1},|\xi_{2}|^{s-2}\xi_{2},\ldots,|\xi_{n}|^{s-2}\xi_{n})+|\xi|^{p-2}\xi.$$
\item Let~$H(x,\xi)=\sqrt[p]{|\xi|_{s}^{p}+\left(\frac{1}{(|x|+1)}+1\right)^{p}|\xi|^{p}}$ ($2\leq s<\infty$) for all~$x\in\Omega$ and~$\xi\in\R^{n}$. Then $H$ satisfies \Cref{ass9}. Suppose that $0\notin\Omega$. Let $\zeta=\infty$ and let~$\{R_{k}\}_{k\in\mathbb{N}}\subseteq (0,\infty)$ be a sequence converging to~$\infty$. Then $\{\mathcal{A}^{(R_{k})}\}_{k\in\N}$ satisfies the conditions in \Cref{limitingd} with the limit operator $$\mathbf{A}(\xi)=|\xi|^{p-s}_{s}(|\xi_{1}|^{s-2}\xi_{1},|\xi_{2}|^{s-2}\xi_{2},\ldots,|\xi_{n}|^{s-2}\xi_{n})+|\xi|^{p-2}\xi.$$ 
\end{itemize}}
\eexample
\blemma
Suppose that~$Q'$ has a Fuchsian-type singularity at~$\zeta$ with respect to a norm $N\in C^{1}(\R^{n}\setminus\{0\})$ on~$\R^{n}$. Then the limiting dilated operator~$Q'_{\{R_{k}\}_{k\in\mathbb{N}}}$ has a Fuchsian-type singularity at both $0$ and $\infty$ with respect to $N$.
\elemma
\bproof
Since~$Q'$ has a Fuchsian-type singularity at~$\zeta$ with respect to $N$, there exist two positive constants $\alpha_{\zeta,\hat{R}}$ and~$\beta_{\zeta,\hat{R}}$ such that for almost all~$x\in\Gamma_{\zeta,\hat{R}}$ and all~$\xi\in\R^{n}$,
   \begin{align*}
\alpha_{\zeta,\hat{R}}|\xi|^{p}\le\mathcal{A}(x,\xi)\cdot\xi\quad\mbox{and}\quad
						|\mathcal{A}(x,\xi)|\le \beta_{\zeta,\hat{R}}\,|\xi|^{{p}-1}.
   \end{align*}
Then for almost all~$x\in\R^{n}\sms\{0\}$ and all~$\xi\in\R^{n}$,
   \begin{align*}
&\alpha_{\zeta,\hat{R}}|\xi|^{p}\le\lim_{k\rightarrow\infty}\mathcal{A}(R_{k}x,\xi)\cdot\xi=\mathbf{A}(x,\xi)\cdot\xi,
\end{align*}
and
\begin{align*}
|\mathbf{A}(x,\xi)|=\lim_{k\rightarrow\infty}\left|\mathcal{A}(R_{k}x,\xi)\right|\le \beta_{\zeta,\hat{R}}\,|\xi|^{{p}-1}.
   \end{align*}
\begin{itemize}
\item Let $p\neq n$. For every $0<a<b<\infty$,
\bals
&\|N(x)^{p-n/q}\mathbf{V}\|_{M^{q}(p;\mathbb{A}_{a,b,N})}\leq b^{p-n/q}\|\mathbf{V}\|_{M^{q}(p;\mathbb{A}_{a,b,N})}\leq b^{p-n/q}\liminf_{k\rightarrow\infty}\|V^{(R_{k})}\|_{M^{q}(p;\mathbb{A}_{a,b,N})}\\
&\leq b^{p-n/q}\liminf_{k\rightarrow\infty}\|R_{k}^{p-n/q}V\|_{M^{q}(p;\mathbb{A}_{aR_{k},bR_{k},N})}\\
&\leq\left(\frac{b}{a}\right)^{p-n/q}\liminf_{k\rightarrow\infty}\|N(x)^{p-n/q}V\|_{M^{q}(p;\mathbb{A}_{aR_{k},bR_{k},N})}.
\eals
Then the following function of $R$,
\bals
&\Lambda_{\mathbf{V}}(R)\triangleq\|N(x)^{p-n/q}\mathbf{V}\|_{M^{q}(p;\mathbb{A}_{a_{N}R,b_{N}R,N})}\\
&\leq\left(\frac{b_{N}}{a_{N}}\right)^{p-n/q}\liminf_{k\rightarrow\infty}\|N(x)^{p-n/q}V\|_{M^{q}(p;\mathbb{A}_{a_{N}RR_{k},b_{N}RR_{k},N})}
\eals
is bounded on $(0,\infty)$.
\item Let $p=n$. For every $0<a<b<\infty$,
\bals
&\|\mathbf{V}\|_{M^{q}(n;\mathbb{A}_{a,b,N})}\leq\liminf_{k\rightarrow\infty}\|V^{(R_{k})}\|_{M^{q}(n;\mathbb{A}_{a,b,N})}\leq\liminf_{k\rightarrow\infty}\|V\|_{M^{q}(n;\mathbb{A}_{aR_{k},bR_{k},N})}.
\eals
Then the following function of $R$,
\bals
&\Lambda_{\mathbf{V}}(R)\triangleq\|\mathbf{V}\|_{M^{q}(n;\mathbb{A}_{a_{N}R,b_{N}R,N})}\leq\liminf_{k\rightarrow\infty}\|V\|_{M^{q}(n;\mathbb{A}_{a_{N}RR_{k},b_{N}RR_{k},N})}
\eals
is bounded on $(0,\infty)$.
\end{itemize}
It follows that the limiting dilated operator $$Q'_{\{R_{k}\}_{k\in\mathbb{N}}}[u]=-\dive\mathbf{A}(x,\nabla u)+\mathbf{V}u^{<p-1>}$$ has a Fuchsian-type singularity at both $0$ and $\infty$ with respect to $N$.
\eproof
The following theorem provides a criterion for the regularity of $\zeta$ with respect to~$Q'$.
\btheorem\label{ldregular}
Suppose that~$Q$ is nonegative in~$\Gamma_{\zeta,\hat{R}}$, that~$V\in\mathbb{M}^{q}_{\loc}(p;\Gamma_{\zeta,\hat{R}})$, that \Cref{C1alpha} holds in~$\Gamma_{\zeta,\hat{R}}$, that~$q>n$ for~$p<2$, and that~$Q'$ has a Fuchsian-type singularity at~$\zeta$ with respect to a norm $N\in C^{1}(\R^{n}\setminus\{0\})$ on~$\R^{n}$. If~$0$ or~$\infty$ is regular with respect to the limiting dilated operator~$Q'_{\{R_{k}\}_{k\in\mathbb{N}}}$, then~$\zeta$ is regular with respect to~$Q'$.
\etheorem
\bproof
The proof is similar to that of \cite[Proposition 3.18]{Giri1}. Let~$u,v\in\mathcal{G}_{\zeta}$. Then $u,v\in\mathcal{G}_{\zeta,\hat{R}}$ up to a smaller (if $\zeta=0$) or larger (if $\zeta=\infty$) $\hat{R}$. By \Cref{monolem}~(i), both $m=\lim_{R\rightarrow|\zeta|}m_{R,N}$ and $M=\lim_{R\rightarrow|\zeta|}M_{R,N}$ exist. 
\begin{itemize}
\item Suppose that~$M=\infty$ or~$m=0$. By \Cref{uhi},~$m=\infty$ or~$M=0$, respectively. Therefore,~$m=M$.
\item Suppose that~$m>0$ and~$M<\infty$. Take~$\tilde{x}\in\R^{n}\setminus\{0\}$. For all sufficiently large~$k\in\N$,~$R_{k}\tilde{x}\in\Gamma_{\gz,\hat{R}}$. For all sufficiently large~$k\in\N$ and all $x\in\Gg_{\gz,\hat{R}}/R_{k}$, we define
$$u_{k}(x)\triangleq\frac{u(R_{k}x)}{u(R_{k}\tilde{x})}\quad\mbox{and}\quad v_{k}(x)\triangleq\frac{v(R_{k}x)}{u(R_{k}\tilde{x})}.$$ For all sufficiently large~$k\in\N$, in~$\Gg_{\gz,\hat{R}}/R_{k}$, both~$u_{k}$ and~$v_{k}$ are positive solutions of the equation
$$-\dive\mathcal{A}^{(R_{k})}(x,\nabla w)+V^{(R_{k})}w^{<p-1>}=0.$$ For all sufficiently large~$k\in\N$,~$u_{k}(\tilde{x})=1$. Since~$m>0$ and~$M<\infty$, for all sufficiently large~$k\in\N$,~$v_{k}(\tilde{x})\asymp 1$. By \Cref{hcpg}, there exists a strictly increasing sequence~$\{k_{j}\}_{j\in\N}$ of positive integers such that~$\{u_{k_{j}}\}_{j\in\N}$ and~$\{v_{k_{j}}\}_{j\in\N}$ converges locally uniformly in~$\R^{n}\sms\{0\}$ to two positive solutions~$u_{\infty}$ and~$v_{\infty}$ of the limiting dilated equation
\begin{align*}
Q'_{\{R_{k}\}_{k\in\mathbb{N}}}[u]=-\dive\mathbf{A}(x,\nabla u)+\mathbf{V}u^{<p-1>}=0,
\end{align*}
respectively. It is elementary to check that~$\{u_{k_{j}}/v_{k_{j}}\}_{j\in\N}$ converges locally uniformly in~$\R^{n}\sms\{0\}$ to~$u_{\infty}/v_{\infty}$. Then for every~$R>0$,
\bals
\max_{S_{R,N}}\frac{u_{\infty}}{v_{\infty}}=\lim_{j\rightarrow\infty}\max_{S_{R,N}}\frac{u_{k_{j}}}{v_{k_{j}}}&=\lim_{j\rightarrow\infty}\max_{x\in S_{R,N}}\frac{u(R_{k_{j}}x)}{v(R_{k_{j}}x)}\\
&=\lim_{j\rightarrow\infty}\max_{x\in S_{RR_{k_{j}},N}}\frac{u(x)}{v(x)}=\lim_{j\rightarrow\infty}M_{RR_{k_{j}},N}=M,
\eals
and
\bals
\min_{S_{R,N}}\frac{u_{\infty}}{v_{\infty}}=\lim_{j\rightarrow\infty}\min_{S_{R,N}}\frac{u_{k_{j}}}{v_{k_{j}}}&=\lim_{j\rightarrow\infty}\min_{x\in S_{R,N}}\frac{u(R_{k_{j}}x)}{v(R_{k_{j}}x)}\\
&=\lim_{j\rightarrow\infty}\min_{x\in S_{RR_{k_{j}},N}}\frac{u(x)}{v(x)}=\lim_{j\rightarrow\infty}m_{RR_{k_{j}},N}=m.
\eals
Since~$0$ or~$\infty$ is regular with respect to the limiting dilated operator~$Q'_{\{R_{k}\}_{k\in\mathbb{N}}}$, $\lim_{x\rightarrow 0}u_{\infty}(x)/v_{\infty}(x)$ or $\lim_{x\rightarrow\infty}u_{\infty}(x)/v_{\infty}(x)$ exists, respectively. Therefore,~$m=M$.
\end{itemize}
Consider any sequence $\{x_{k}\}_{k\in\N}\subseteq\Gamma_{\zeta,\hat{R}}$ converging to $\zeta$. Then $\lim_{k\rightarrow\infty}N(x_{k})=|\zeta|$. Therefore,
$$\limsup_{k\rightarrow\infty}\frac{u(x_{k})}{v(x_{k})}\leq\lim_{k\rightarrow\infty}M_{N(x_{k}),N}=M=m=\lim_{k\rightarrow\infty}m_{N(x_{k}),N}\leq\liminf_{k\rightarrow\infty}\frac{u(x_{k})}{v(x_{k})}.$$
It follows that $\lim_{k\rightarrow\infty}u(x_{k})/v(x_{k})=M=m$.
Then~$\lim_{x\rightarrow\zeta}u(x)/v(x)$ exists, which means that~$\zeta$ is regular with respect to~$Q'$.
\eproof
Now we define the weak Fuchsian-type singularity.
\bdefinition
\emph{Suppose that~$Q'$ has a Fuchsian-type singularity at~$\zeta_{1}=\zeta$. Let~$m\in\N$. Let~$\{R_{k,1}\}_{k\in\N}$ be as in \Cref{limitingd}. For every~$2\leq j\leq m$, let~$\{R_{k,j}\}_{k\in\N}$ converge to~$\zeta_{j}\in\{0,\infty\}$ and satisfy all the other conditions in \Cref{limitingd} with respect to~$\mathbf{A}_{\{R_{k,j-1}\}_{k\in\N}}$ and~$\mathbf{V}_{\{R_{k,j-1}\}_{k\in\N}}$. Suppose that for every $1\leq j\leq m$, there exists a norm family $\mathbf{H}_{\{R_{k,j}\}_{k\in\N}}$ satisfying \Cref{ass9} such that $\mathbf{A}_{\{R_{k,j}\}_{k\in\N}}(x,\xi)=\nabla_{\xi}\left(\mathbf{H}_{\{R_{k,j}\}_{k\in\N}}(x,\xi)^{p}/p\right)$ for almost all $x\in\R^{n}\sms\{0\}$ and all $\xi\in\R^{n}$. The operator~$Q'$ has a \emph{weak Fuchsian-type singularity} at~$\zeta$ if in~$\R^{n}\sms\{0\}$, $\mathbf{H}_{\{R_{k,m}\}_{k\in\N}}$ is a fixed norm and~$\mathbf{V}_{\{R_{k,m}\}_{k\in\N}}=0$.}
\edefinition
With these preparations at our disposal, we can prove the following positive Liouville-type theorem for $p\leq n$.
\btheorem[{Positive Liouville-type theorem}]\label{Liouville}
Let~$p\leq n$ and let $N\in C^{1}(\R^{n}\sms\{0\})$ be a norm on $\R^{n}$. Suppose that~$Q'$ has a Fuchsian-type singularity at $\zeta$ with respect to $N$ and a weak Fuchsian-type singularity at~$\zeta$. For every~$j=1,2,\cdots,m$ and some~$\hat{R}_{j}\in(0,\infty)$, where $\hat{R}_{1}=\hat{R}$, suppose that \Cref{C1alpha} holds in~$\Gamma_{\zeta_{j},\hat{R}_{j}}$ with respect to $\mathbf{A}_{\{R_{k,j-1}\}_{k\in\N}}$, that~$Q_{\mathbf{A}_{\{R_{k,j-1}\}_{k\in\N}},\mathbf{V}_{\{R_{k,j-1}\}_{k\in\N}}}$ is nonegative in~$\Gamma_{\zeta_{j},\hat{R}_{j}}$, where $\mathbf{H}_{\{R_{k,0}\}_{k\in\N}}=H$ and $\mathbf{V}_{\{R_{k,0}\}_{k\in\N}}=V$, that $\mathbf{V}_{\{R_{k,j-1}\}_{k\in\N}}\in\mathbb{M}^{q}_{\loc}(p;\Gamma_{\zeta_{j},\hat{R}_{j}})$, and that~$q>n$ for~$p<2$. Suppose that~$\mathbf{H}_{\{R_{k,m}\}_{k\in\N}}$ satisfies \Cref{C1alphan}. Then~$\zeta$ is regular with respect to~$Q'$ and~$\mathbf{M}_{\zeta,\hat{R}}$ has a unique function up to a positive multiplicative constant.
\etheorem
\bremark
\emph{This theorem still holds if we assume Fuchsian-type singularities for $Q'$ and $Q'_{\mathbf{A}_{\{R_{k,j-1}\}_{k\in\N}},\mathbf{V}_{\{R_{k,j-1}\}_{k\in\N}}}$ only on an essential subset of~$\Gamma_{\zeta,\hat{R}}$ at~$\zeta$ and only on an essential subset of~$\Gamma_{\zeta_{j},\hat{R}_{j}}$ at~$\zeta_{j}$ for every~$j=1,2,\cdots,m$ with respect to $N$, respectively.}
\eremark
\bproof[Proof of \Cref{Liouville}:]
The proof is similar to that of \cite[Theorem 5.4]{Giri1}. By \Cref{pncoro},~$0$ is regular with respect to $Q'_{\mathbf{A}_{\{R_{k,m}\}_{k\in\N}},0}$. By \Cref{ldregular},~$\zeta_{m}$ is regular with respect to $Q'_{\mathbf{A}_{\{R_{k,m-1}\}_{k\in\N}},\mathbf{V}_{\{R_{k,m-1}\}_{k\in\N}}}$. By induction, we deduce that~$\zeta_{1}=\zeta$ is regular with respect to $Q'$. Then \Cref{uniquemini} yields the second claim.
\eproof
\appendix
\section*{Acknowledgments}
This work is part of the author’s ongoing Ph.D. thesis of the Technion - Israel Institute of Technology written under joint supervision of Professors Yehuda Pinchover and Matthias Keller. The author thanks Professor Pinchover for his authoritative guidance. The author also acknowledges the financial support provided by the Technion and the Israel Science Foundation (Grant No. 637/19) founded by the Israel Academy of Sciences and Humanities.
{}

\end{document}